\documentstyle[11pt]{article}
\input{amssymb.sty}
\input{psfig.sty}
\normalbaselines
\begin{document}
\bibliographystyle{unsrt}
\newtheorem{thm}{{\sc Theorem}}[section]
\newtheorem{lma}{{\sc Lemma}}[section]
\newtheorem{rmk}{{\sc Remark}}[section]
\newtheorem{prp}{{\sc Proposition}}[section]
\newtheorem{cor}{{\sc Corollary}}
\newtheorem{lemma}{{Lemma}}
\newtheorem{theorem}{{Theorem}}

\newfont{\msbm}{msbm10 scaled\magstep1}
\newfont{\eusm}{eusm10 scaled\magstep1}

\def\vv{{\bf v}}
\def\vr{{\bf r}}
\def\vp{{\bf p}}
\def\spec{{\rm spec}}
\def\Hom{{\rm Hom}}
\def\hl{{\rm Hilb}}
\def\NZ{{\mathbb N}}
\def\ZZ{{\mathbb Z}}
\def\RZ{{\mathbb R}}
\def\CZ{{\mathbb C}}
\def\PZ{{\mathbb P}}
\def\QZ{{\mathbb Q}}

\def\eM{{J_0^{\dagger}}}
\def\Io{{I(\widehat{o})}}
\def\bas{{\mbox{\eusm B}}}
\def\wt{{\rm wt}}
\def\RT{{\goth R}}

\def\bea*{\begin{eqnarray*}}
\def\eea*{\end{eqnarray*}}
\def\ba{\begin{array}}
\def\ea{\end{array}}
\count1=1
\def\be{\ifnum \count1=0 $$ \else \begin{equation}\fi}
\def\ee{\ifnum\count1=0 $$ \else \end{equation}\fi}
\def\ele(#1){\ifnum\count1=0 \eqno({\bf #1}) $$ \else \label{#1}\end{equation}\fi}
\def\req(#1){\ifnum\count1=0 {\bf #1}\else \ref{#1}\fi}
\def\bea(#1){\ifnum \count1=0   $$ \begin{array}{#1}
\else \begin{equation} \begin{array}{#1} \fi}
\def\eea{\ifnum \count1=0 \end{array} $$
\else  \end{array}\end{equation}\fi}
\def\elea(#1){\ifnum \count1=0 \end{array}\label{#1}\eqno({\bf #1}) $$
\else\end{array}\label{#1}\end{equation}\fi}
\def\cit(#1){
\ifnum\count1=0 {\bf #1} \cite{#1} \else 
\cite{#1}\fi}
\def\bibit(#1){\ifnum\count1=0 \bibitem{#1} [#1    ] \else \bibitem{#1}\fi}
\def\ds{\displaystyle}
\def\hb{\hfill\break}
\def\comment#1{\hb {***** {\em #1} *****}\hb }

\def\goth{\mathfrak}
\newtheorem{pf}{Proof}
\renewcommand{\thepf}{}
\vbox{\vspace{38mm}}
\begin{center}
{\LARGE \bf  On Hypersurface Quotient Singularities of Dimension 4 }\\[5mm]
Li Chiang 
\\{\it Institute of Mathematics \\ Academia Sinica \\ 
Taipei , Taiwan \\
(e-mail:chiangl@gate.sinica.edu.tw)}
\\[5mm]
Shi-shyr Roan
\\{\it Institute of Mathematics \\ Academia
Sinica \\  Taipei , Taiwan \\ (e-mail:
maroan@ccvax.sinica.edu.tw)} \\[5mm]
\end{center}

\begin{abstract} 
We consider geometrical problems on 
Gorenstein hypersurface orbifolds of dimension
$n \geq 4$ through the theory of Hilbert scheme of
group orbits. For a linear special group $G$
acting on
$\CZ^n$, we study the $G$-Hilbert
scheme, $\hl^G(\CZ^n)$, and crepant
resolutions of $\CZ^n/G$ for $G$=the $A$-type
abelian group $ A_r(n)$. For
$n=4$, we obtain the explicit structure of 
$\hl^{A_r(4)}(\CZ^4)$. The crepant
resolutions of $\CZ^4/A_r(4)$ are constructed through
their relation with $\hl^{A_r(4)}(\CZ^4)$,
and the connections between these crepant
resolutions are found by the ``flop" procedure of 
4-folds. We also make some primitive discussion
on $\hl^G(\CZ^n)$ for 
the $G$= alternating group ${\goth A}_{n+1}$ of
degree $n+1$ with the standard representation on
$\CZ^n$;  the detailed structure of  $\hl^{{\goth A}_4}(\CZ^3)$ is explicitly constructed.
\par \vspace{5mm} \noindent

\end{abstract}

\vfill
\eject

\section{Introduction}
The purpose of this paper is to study some 
geometrical problems of certain
 Gorenstein hypersurface 
orbifolds of dimension 4.
The main focus will be on the
structure of the newly developed concept of Hilbert
scheme of group orbits and its connection with
crepant resolutions of the orbifold. For a finite
subgroup $G$ in ${\rm SL}_n(\CZ)$, the 
$G$-Hilbert scheme,
$\hl^G(\CZ^n)$, was first introduced by 
Nakamura et al
\cite{GNS, IN1, IN2, N1}; one primary goal aims 
to provide a conceptual understanding of 
crepant resolutions of $\CZ^n/G$ for $n=3$, whose
solution was previously known by a 
computational method, relying heavily on  
Miller-Blichfeldt-Dickson classification of
finite groups in ${\rm SL}_3(\CZ)$ \cite{MBD} and
the invariant theory of two simple groups,
$I_{60}$ (icosahedral group), 
$H_{168}$ (Klein group) \cite{Kl} ( for the 
existence of crepant resolutions, see 
\cite{Rtop} and references therein). For $n=2$, 
$\hl^G(\CZ^2)$ is 
the minimal resolution of $\CZ^2/G$, hence it has
the trivial canonical bundle \cite{IN1, IN2, N1}. For $n=3$,
it was expected that $\hl^G(\CZ^3)$ is a
crepant resolution of $\CZ^3/G$. Recently the
affirmative answer was obtained, in \cite{INj, N} for the
abelian group
$G$  by techniques in toric
geometry,  and in \cite{BKR} for a general group $G$
by derived category methods bypassing the
geometrical analysis of $G$-Hilbert scheme.
With these successful results in dimension 3, a
question naturally arises on the possible role of
$G$-Hilbert scheme on crepant resolution problems of
orbifolds in dimension $n \geq 4$.
For $n \geq 4$, it is a well-known fact that 
$\CZ^n/G$ might have no 
crepant resolutions at all, even
for a cyclic group $G$ and $n=4$, (for a
selection of examples, see e.g., \cite{R97}). 
To avoid many such 
complicated exceptional cases,  we
will restrict our attention only to those with
hypersurface singularities. 
In this paper, we will 
address certain problems on two specific types of
hypersurface Gorenstein quotient singularity,
$\CZ^n/G$, of dimension $n$; one is the abelian
group $G= A_r(n)$ defined in (\ref{Arn}) in the
main body of the paper, the other group $G$ is  the
alternating group ${\goth A}_{n+1}$ of degree
$n+1$ acting on
$\CZ^n$ through the standard representation. In
the case
$G= A_r(n)$, $\hl^G(\CZ^n)$ is a toric
variety, hence the methods for toric
geometry provide an effectively
tool to study its structure.  
For $n=4$, we will give a detailed derivation
of the  smooth toric structure of ${\rm
Hilb}^{A_r(4)}(\CZ^4)$, and construct
the crepant toric resolutions of
$\CZ^4/A_r(4)$ by blowing-down the canonical
divisors of $\hl^{A_r(4)}(\CZ^4)$; in due
course the ``flop" of 4-folds 
naturally arises in the process, (see Theorems
\ref{th:A1(4)}, \ref{th:flip}, \ref{th:A(4)} in
the main text of this paper, whose statements 
were previously announced in \cite{CR}). We
would expect the concept appeared in the proof of
these theorems will inspire certain clue to other
cases, not only the
$A_r(n)$-type groups, but also for the non-abelian
groups ${\goth A}_{n+1}$ ( which are simple groups
for $n \geq 4$). The group ${\goth A}_4$
is a solvable group of order 12, also called the
ternary trihedral group. The crepant
resolution of
$\CZ^3/{\goth A}_4$ was explicitly constructed in
\cite{BM}, and the structure ${\rm
Hilb}^{{\goth A}_4}(\CZ^3)$ over the origin orbit
of $\CZ^3/{\goth A}_4$ was analyzed in detail in
\cite{GNS}. Through the
representation theory of ${\goth A}_4$, we will give the direct verification
that ${\rm
Hilb}^{{\goth A}_4}(\CZ^3)$ is smooth and a crepant resolution of 
$\CZ^3/{\goth A}_4$.  Though
the conclusion is known by the general result
in \cite{BKR} using qualitative arguments, the object of our detailed analysis
aims to reveal that there exist certain common features in determining the structures of
$G$-Hilbert schemes for certain
abelian and non-abelian groups $G$ by the computational methods, in hope that
 the approach could possibly
be applied to higher
dimensional cases. With this in mind,  we will in
this paper restrict our attention only to  the
case ${\goth A}_4$, leave possible
generalizations, applications or implications to
future work.

This paper is
organized as follows.  In \S2, we will
summarize the main features of  $G$-Hilbert
scheme of group orbits, and results in toric
geometry for the need of later discussion. We
will also define the group $G$ which we will consider in
this paper. The next two sections will be devoted
to the discussion of the structure of
$\hl^G (\CZ^4)$ and crepant resolutions of
$\CZ^4/G$ for $G=
A_r(4)$. For the simpler terminology to express
the ideas, also for the description of
geometry of flop of 4-folds,  we will consider only the case
$A_1(4)$ in
\S3 to discuss the structure of ${\rm
Hilb}^{A_1(4)}(\CZ^4)$. The flop relation between
crepant resolutions of
$\CZ^4/A_1(4)$ will be examined in detail through
$\hl^{A_1(4)}(\CZ^4)$. In \S4, we will
derive the solution of the corresponding problems
for
$G= A_r(4)$ for a general positive integer $r$,
with much more complicated techniques but a
method much in tune with the previous section. In
\S5, we consider the case
$G= {\goth A}_{n+1}$ acting on
$\CZ^n$ through the standard representation for $n=3$.
By employing the structure of the fiber 
in $\hl^{{\goth A}_4}(\CZ^3)$ 
over the origin orbit of $\CZ^3/{\goth A}_4$
described in
\cite{GNS}, we give an explicit
construction of the smooth and crepant structure 
of $\hl^{{\goth A}_4}(\CZ^3)$ using finite
group representation theory, along a line similar to
 the previous two sections in a certain sense.
Finally we give the conclusion remarks in \S6.

{\bf Notations. } To present 
our work, we prepare some notations. In this
paper, by an orbifold we shall always mean the
orbit space of a smooth complex 
manifold acted on by a finite
group. Throughout the paper, $G$ will always
denote a finite group unless otherwise stated. We
denote 
\begin{eqnarray*}
{\rm Irr}(G) &= & \{ \rho : G \longrightarrow
{\rm GL}( V_{\rho}) \ \  {\rm an
\ irreducible
\ representation \ of \ } G \} . 
\end{eqnarray*}
The trivial representation of $G$ will be denoted
by ${\bf 1}$. For a $G$-module $W$,
i.e., a $G$-linear representation 
space 
$W$, one has the canonical irreducible
decomposition: $
W = \bigoplus_{\rho \in {\rm Irr}(G)} W_{\rho}$, 
where $W_{\rho}$ is a $G$-submodule of
$W$, isomorphic to $V_{\rho} \otimes
W_{\rho}^0$ for some  trivial $G$-module
$W_{\rho}^0$. 
For an analytic variety $X$, we shall not distinguish the notions of 
vector bundle and locally free ${\cal O}_X$-sheaf over $X$.

\section{G-Hilbert Scheme, Toric Geometry}
In this section, we brief review some basic facts
on $\hl^G(\CZ^n)$ ( the Hilbert scheme of
$G$-orbits) and toric geometry necessary for
later use, then specify the groups $G$
for the discussion of the rest sections of this paper.

First, we will always assume 
$G$ to be a finite subgroup of
${\rm SL}_n(\CZ)$. Denote 
$S_G:=\CZ^n/G$ with the canonical projection,
$
\pi_G:\CZ^n\rightarrow S_G $, 
and $o:=\pi_G(\vec{0})$.
As $G$ acts on
$\CZ^n$ freely outside a finite collection of
linear subspaces with codimension $\ge 2$, $S_G$
is an orbifold with non-empty singular set
${\rm Sing}(S_G)$ of codimension $\ge 2$. In fact,
the element $o$ is a
singular point of $S_G$.
By a variety $X$ birational over $S_G$, we will always mean
a proper birational morphism  $\sigma$ from 
$X$ to $S_G$   which is  biregular 
between
$X
\setminus \sigma^{-1}({\rm Sing}(S_G))$ and $ S_G
\setminus {\rm Sing}(S_G) $,
\be
\sigma : X \longrightarrow  S_G \ .
\ele(sigma)
One can form the commutative diagram 
via the birational morphism $\sigma$,
\bea(cll)
X \times_{S_G} \CZ^n & \longrightarrow & \CZ^n \\ 
\ \ \downarrow \pi & & \downarrow \pi_G  \\
X  &\stackrel{\sigma}{\longrightarrow} & S_G  \ .
\elea(Xdia)
Denote ${\cal F}_X$ the coherent 
${\cal O}_X$-sheaf over $X$ obtained by 
the push-forward of the structure sheaf 
of $X \times_{S_G} \CZ^n$, $
{\cal F}_X : = \pi_* {\cal O}_{X \times_{S_G}\CZ^n}$.
For two varieties $X,
X^\prime$  birational over $S_G$ with
the commutative diagram,
$$
\begin{array}{lll}
\ \ X  &\stackrel{\sigma}{\longrightarrow}
& S_G \\ 
\mu \downarrow & &  || \\
\ \ X^\prime
&\stackrel{\sigma^\prime}{\longrightarrow} & S_G  \ ,
\end{array}
$$
one has a canonical morphism, 
$ \mu^*{\cal F}_{X^\prime} \longrightarrow {\cal
F}_X$.  In particular, the morphism
(\req(sigma)) gives rise to the ${\cal O}_X$-morphism,
$$
\sigma^* ( {\pi_G}_*{\cal O}_{\CZ^n} )
\longrightarrow {\cal F}_X  \ .
$$
Furthermore, all the above morphisms  are
compatible with the natural $G$-structure of
${\cal F}_X$ induced from the
$G$-action on
$\CZ^n$ via (\req(Xdia)). Then ${\cal
F}_X$ has the
canonical $G$-decomposition of coherent ${\cal O}_X$-submodules:
$
{\cal F}_X = \bigoplus_{\rho \in {\rm Irr}(G)}
({\cal F}_X)_\rho $,
where $({\cal F}_X)_\rho$ is the $\rho$-factor of 
${\cal F}_X$. The geometrical fibers of ${\cal
F}_X$ and $ ({\cal F}_X)_\rho$ over $x \in X$ are defined by
$
{\cal F}_{X, x} = k(x) \bigotimes_{{\cal O}_X} {\cal F}_X$, $ 
({\cal F}_X)_{\rho, x} = k(x) \bigotimes_{{\cal
O}_X} ({\cal F}_X)_\rho $, 
where $k (x) ( := {\cal O}_{X, x}/{\cal M}_{x})$
is the residue field  at $x$. Over 
$X \setminus \sigma^{-1}({\rm Sing}(S_G))$, 
${\cal F}_X$ 
is a vector bundle of rank $|G|$  with the
regular $G$-representation on each geometric
fiber. Hence $({\cal F}_X)_\rho$ is a vector
bundle over 
$X \setminus \sigma^{-1}({\rm Sing}(S_G))$ with
the  rank equal to the dimension of $ V_\rho$. For $x \in X$, there
exists a
$G$-invariant ideal $I(x)$ in $\CZ[Z] (:=
\CZ[Z_1, \cdots, Z_n ])$ such that 
the following relation holds,
\begin{eqnarray}
{\cal F}_{X, x} = k (x)\bigotimes_{{\cal O}_{S_G} }
{\cal O}_{\CZ^n} 
(x) \simeq \CZ[Z]/I(x) \ .
\label{I(x)} 
\end{eqnarray}
We have $({\cal F}_X)_{\rho, x}
\simeq (\CZ[Z]/I(x))_\rho$. 
The vector spaces
$\CZ[Z]/I(x)$ form a family of finite
dimensional $G$-modules parametrized by
$x \in X$. For $ x \not\in \sigma^{-1}({\rm
Sing}(S_G))$,  $\CZ[Z]/I(x)$ is a 
regular $G$-module.
In particular, for  $X=S_G$ in (\ref{I(x)}) and $s \in
S_G$, the 
$G$-invariant ideal $I(s)$ of $\CZ[Z]$  is generated by the $G$-invariant
polynomials vanishing at $\sigma^{-1}(s)$.
Let $\widetilde{I}(s)$  be the ideal of $\CZ[Z]$
consisting of all polynomials 
vanishing at
$\sigma^{-1}(s)$. Then $\widetilde{I}(s)$ is an
$G$-invariant ideal with $
\widetilde{I}(s) \supset I(s)$. 
For $s= o$, we have $  
\widetilde{I}(o) = \CZ[Z]_0$ and $ 
I(o) = \CZ[Z]^G_0 \CZ[Z]$, 
where the subscript $0$ indicates the maximal
ideal of polynomials vanishing at the
origin.  
For a variety $X$ birational over $S_G$
via $\sigma$ in (\req(sigma)), one has the
following relations of $G$-invariant ideals of
$\CZ[Z]$:
$$
\widetilde{I}(s) \supset I(x) \supset
I(s) \ , \ \ x \in X \ , \ s = \sigma(x) \ . 
$$ 
For $x \in X$, 
there exists a direct sum decomposition of
$\CZ[Z]$ as $G$-modules,
$$
\CZ[Z] = I(x)^{\bot} \oplus I(x) \ .
$$
Here $I(x)^{\bot}$ is a finite dimensional
$G$-module isomorphic to
$\CZ[Z]/I(x)$. Similarly,  we have $G$-module decompositions for 
$s=\sigma(x) \in S_G$,
$$
\CZ[Z] = I(s)^{\bot} \oplus
I(s) \ ,
\ \ \CZ[Z] = \widetilde{I}(s)^{\bot} \oplus
\widetilde{I}(s)
$$ 
so that the relations,
$\widetilde{I}(s)^{\bot} \subset I(x)^{\bot}
\subset I(s)^{\bot}$, hold. Note that the above
finite dimensional $G$-modules with superscript
$\bot$ are not unique in $\CZ[Z]$ because there is a choice involved,
nonetheless we could choose them such that this inclusions are fulfilled. One has the canonical
$G$-decomposition of $I(x)^{\bot}$: $
I(x)^{\bot} = \bigoplus_{\rho \in {\rm Irr}(G)}
I(x)^{\bot}_\rho$, 
where the factor $I(x)^{\bot}_\rho $ is isomorphic
to a positive finite sum of copies of
$V_{\rho}$.

Now we consider the varieties       
$X$ birational over $S_G$ such that ${\cal F}_X$ is a
vector bundle. Among all such $X$,   there
exists a  minimal object, called the $G$-Hilbert
scheme in
\cite{IN1, IN2, N1, N},
\be
\sigma_\hl : \hl^G(\CZ^n)
\longrightarrow  S_G \  .
\ele(sHilb) 
By the definition of $\hl^G(\CZ^n)$, an
element (i.e. closed point) $p$ of $\hl^G(\CZ^n)$ is described by a
$G$-invariant ideal $I$($=I(p)$) of $\CZ[Z]$ of
colength $|G|$,  and the fiber of the vector
bundle ${\cal F}_{\hl^G(\CZ^n)}$ over $p$
can be identified with the regular $G$-module $
\CZ[Z]/I(p)$. For simplicity of notations, we shall also make the identification of the element $p$ with its associated ideal $I$, and write $I\in \hl^G(\CZ^n)$ in what follows if no confusion arises.
For any other $X$, the map (\req(sigma)) can be factored through 
a birational morphism $\lambda$ from $X$ onto
$\hl^G(\CZ^n)$ via $\sigma_\hl$,
$$
\begin{array}{lll}
\ \ X  &\stackrel{\sigma}{\longrightarrow}
& S_G \\ 
\lambda \downarrow & &  || \\
\ \ \hl^G(\CZ^n)
&\stackrel{\sigma_\hl}{\longrightarrow} &
S_G  \ .
\end{array}
$$
In fact, the ideal $I (x)$ of
(\ref{I(x)}) is a  colength $|G|$ ideal in
$\CZ[Z]$,  by which the map 
$\lambda$ is defined. We
will denote $X_G$ the normalization of ${\rm
Hilb}^G(\CZ^n)$, which is a normal variety
over $S_G$ with the birational morphism from $ X_G$ onto $S_G$.
As  every biregular automorphism of
$S_G$ can always be lifted to one of 
$\hl^G(\CZ^n)$, hence also to $X_G$,
one has the following result.
\begin{lma} \label{lem:AutSG}
Let ${\rm Aut}(S_G)$ be the group of biregular automorphisms 
of $S_G$. Then $\hl^G(\CZ^n), X_G$ are 
${\rm Aut}(S_G)$-varieties over $S_G$ via ${\rm
Aut}(S_G)$-morphisms. As a consequence,
$X_G$ is a toric variety for an abelian group $G$.
\end{lma}

Now we are going to summarize some basic facts in
toric geometry for the later discussion when the group $G$ is abelian,
( for details, see e.g., \cite{D, M, O}) . In this
case, we consider
$G$ as a subgroup of the diagonal group $T_0$ of
${\rm GL}_n(\CZ^n)$ with the identification
$T_0 =
\CZ^{* n}$. Regard 
$\CZ^n$ as the partial compactification of $T_0$,
\begin{eqnarray*}
G \subset T_0 
\subset \CZ^n \ .
\end{eqnarray*}
Let $T$ be the torus $T_0/G$ and consider $S_G \
(= \CZ^n/G)$ as a $T$-space, 
$$
T := T_0/G \ , \ \ \ T \subset S_G \ .
$$
The combinatorial data of toric varieties are
constructed from the lattices of 1-parameter
subgroups and characters of tori $T, T_0$,
\begin{eqnarray*}
N (:= {\rm Hom}(\CZ^*, T )) &\supset& N_0 (:= {\rm
Hom}(\CZ^*, T_0 )) \ , \\
M (:= {\rm Hom}( T, \CZ^* )) &\subset& M_0 (:=
{\rm Hom}( T_0, \CZ^* )) \ .
\end{eqnarray*}
For convenience, $N_0, N$ will
be identified with the following lattices in
$\RZ^n$ in this paper. Denote by $\{ e^i \}_{i=1}^n$ the 
standard basis of $\RZ^n$, and define the map
$ 
{\rm exp}: {\RZ}^n \longrightarrow T_0 $ by $ 
r (=\sum_{i=1}^n r_ie^i ) \mapsto {\rm
exp}(r) := \sum_{i} e^{2\pi
\sqrt{-1}r_i}e^i$. 
The lattices $N, N_0$ are given by 
$$
N_0  = \ZZ^n (:= {\rm exp}^{-1}(1)) \ , \ \ 
N = {\rm exp}^{-1}( G ) \ , 
$$
and we have $G \simeq N/N_0$.
The lattice $M_0$ dual to $N_0$ is the standard
one in the dual space $\RZ^{n *}$. In what
follows, we shall identify $M_0$ with the group of
monomials in  variables
$Z_1, \ldots, Z_n$ via the correspondence:
$$
I = \sum_{s=1}^n i^se_s \in M_0 \ \ 
\longleftrightarrow \ \ Z^I = \prod_{s=1}^n
Z_s^{i_s} \ .
$$  
The dual lattice $M$ of $N$ is the sublattice of 
$M_0$, consisting of all $G$-invariant
monomials.
Among the varieties 
$X$  birational over the $T$-space $S_G$, we shall
consider only  those 
$X$ with a $T$-structure.
It has been known that these toric varieties $X$
are represented by certain  combinatorial data in
toric geometry. A toric variety over $S_G$ is
described by  a fan
$\Sigma = \{ \sigma_{\alpha} \ | \ \sigma \in I
\}$ with  the first
quadrant of
$\RZ^n$ as its support, i.e., 
a rational convex cone decomposition of the first
quadrant in
$\RZ^n$. Equivalently, these combinatorial data
can also be described by the intersection of the
fan and the standard simplex
$\Delta$ in the first quadrant,
\be
\bigtriangleup : = \{ r \in {\RZ}^n | \sum_{i}
r_i = 1, r_j \geq 0 \ \ \forall \ j \} \ .
\ele(tri)
The corresponding data in $\bigtriangleup $ are
denoted by
$\Lambda  = \{
\bigtriangleup_{\alpha} \ | \ \alpha \in I \} $
with
$
\bigtriangleup_{\alpha}:= \sigma_{\alpha} \cap 
\bigtriangleup$. Then  $\Lambda $
is a polytope decomposition of
$\bigtriangleup$ with vertices in $\bigtriangleup \cap
\QZ^n$. Note that for $\sigma_{\alpha}=
\{ \vec{0} \}$, we have
$\bigtriangleup_{\alpha}= \emptyset$. Such 
$\Lambda $ will be called a
rational polytope decomposition of
$\bigtriangleup$, and we will denote
$X_\Lambda$ the toric variety corresponding to
$\Lambda$. If all vertices of 
$\Lambda$ are in $N$, 
$\Lambda$ is called  an integral polytope
decomposition of
$\bigtriangleup$. For a 
rational polytope decomposition $\Lambda$ of
$\bigtriangleup$, we define $\Lambda(i):= \{
\Delta_{\alpha} \in \Lambda \ | \ {\rm
dim}(\Delta_{\alpha}) = i \}$ for $ -1 \leq i
\leq n-1$, (here ${\rm dim}( \emptyset ):= -1$).
The
$T$-orbits in
$X_\Lambda$ are parametrized by
$\bigsqcup_{i=-1}^{n-1} \Lambda(i)$.  In fact,
for 
$\bigtriangleup_{\alpha} \in \Lambda(i)$, there
associates a $(n-1-i)$-dimensional
$T$-orbit, which will be denoted
by ${\rm orb}(\bigtriangleup_{\alpha})$.  A toric
divisor in $X_\Lambda$ is the closure of an 
$(n-1)$-dimensional orbit, denoted by $D_v =
\overline{{\rm orb}(v)}$ for $v \in 
\Lambda (0)$. The canonical sheaf of 
$X_{\Lambda}$ is expressed by the toric divisors
(see, e.g.
\cite{D, M, O}),
\begin{equation}
\omega_{X_\Lambda} = {\cal O}_{X_\Lambda}(
\sum_{v \in \Lambda(0) } (m_v-1)D_v ) \ , \ 
\label{K}
\end{equation}
where $m_v$ is the  least positive integer 
with  $m_v v \in N$. In particular, $X_\Lambda$ is crepant, i.e. ,  
$\omega_{X_\Lambda}= {\cal O}_{X_\Lambda}$, if and only if
 $\Lambda$ is integral.
On the other hand, the smoothness of  $X_{\Lambda}$
is described by  the decomposition
$\Lambda$ to be a simplicial one with
the multiplicity one property, i.e., for each
$\Lambda_{\alpha} \in
\Delta(n-1)$, the elements $m_v v$ for $v \in
\Lambda_{\alpha} \cap \Lambda(0)$  form a
$\ZZ$-basis of $N$. 
The
following results are known for toric variety
over $S_G$  (see e.g. \cite{R89}):

(1) The Euler number of $X_\Lambda$ is given by
$\chi ( X_\Lambda ) = | \Lambda (n-1) |$.

(2) For a rational polytope
decomposition $\Lambda$ of $\Delta$, any
two of the following three
properties imply the third
one:
\begin{eqnarray} 
X_\Lambda : {\rm non-singular
} , \ \ \ \omega_{X_\Lambda} = {\cal O}_{X_\Lambda} \
, \ \ \ \chi ( X_\Lambda ) = | G | . \label{tori}
\end{eqnarray}

In this paper, we shall consider only two
specific series of hypersurface $n$-orbifold
$S_G$ for $n \geq 2$. 
The first
type can be regarded as a
generalization of the $A$-type Klein surface
singularity, the group
$G$ is defined as follows,
\begin{eqnarray}
A_r(n) := \{ g \in {\rm SL}_n(\CZ) \ | \ g:
{\rm diagonal} \ ,\  g^{r+1}=1  \} \ , \ \
r \geq 1 \ . \label{Arn}
\end{eqnarray}
The 
$A_r(n)$-invariant polynomials in $\CZ[Z] (:=
\CZ[Z_1, \ldots, Z_n])$ are generated by
monomials, $
X := \prod_{i=1}^n Z_i$ and $  Y_j :=
Z_j^{r+1} \ (j=1,
\ldots, n)$. 
Thus $S_{A_r(n)}$ is realized as the 
hypersurface in
$\CZ^{n+1}$,
$$
S_{A_r(n)} : \ \ \ \ \  x^{r+1} = \prod_{j=1}^n 
y_j \ , \ \ \ \ (x, y_1, \ldots, y_n ) \in
\CZ^{n+1} \ .
$$
The ideal $I(o)$ of 
$\CZ[Z]$ for the element $o \in
S_{A_r(n)}$ is given by $
I(o) = \langle Z_1^{r+1}, \ldots, Z_n^{r+1}, Z_1 \cdots
Z_n \rangle $, 
hence
$$
I(o)^{\bot} = \bigoplus \{ \CZ Z^I \ | \ I=(i^1,
\ldots, i^n) \ , \ \ 0 \leq i^j \leq r \ ,
\prod_{j=1}^n i^j=0 
\} \ .
$$
For a nontrivial character $\rho$ of $A_r(n)$,
the dimension of $I(o)^{\bot}_\rho$ is always
greater than one. In fact, one can describe 
an explicit set of monomial generators of 
$I(o)^{\bot}_\rho$. For example, say
$I(o)^{\bot}_\rho$ containing an element $Z^I$ with 
$I=(i^1, \ldots, i^n), i^1=0$ and $ i^s \leq
i^{s+1}$, then $I(o)^{\bot}_\rho$ is generated by
$Z^K$s with 
$K=(k^1,  \ldots, k^n)$ given by   
\be
k^s= \left\{ \begin{array}{ll} 
r+1-i^j + i^s ,  & {\rm if} \  i^s <i^j \ , \\
i^s - i^j , & {\rm otherwise} \ ,
\end{array} \right. 
\ele(Ioorg)
here  $j$ runs through $1$ to $ n$. Note that
some of the above $n$-tuples $K$ might coincide. In
particular for
$r=1$, the dimension of 
$I(o)^{\bot}_\rho$ is equal to $2$ for $\rho 
\neq {\bf 1}$, with a basis consisting of $Z^I, Z^{I^\prime}$
whose indices satisfy the relations, 
$0 \leq i^s,i^{s^\prime}\leq 1 , i^s+i^{s^\prime}=1$ for $1 \leq s \leq n$. 

The second type of group $G$ is the alternating
group ${\goth A}_{n+1}$ (of degree
$n+1$) acting on
$\CZ^n$ through the standard representation. The representation is induced from the
linear action of the symmetric group
${\goth S}_{n+1}$ on
$\CZ^{n+1}$ by permuting the coordinate indices,
then restricting on the subspace  
\be
V = \{ (\tilde{z}_1, \ldots, \tilde{z}_{n+1}) \in
\CZ^{n+1}
\ |
\ 
\sum_{j=1}^{n+1} \tilde{z}_j = 0 \} \ \simeq \
\CZ^n \ \ .
\ele(V)
We denote $\CZ[ \tilde{Z}] ( := \CZ[\tilde{Z}_1,
\cdots ,
\tilde{Z}_{n+1} ])$ the coordinate ring of the
affine $(n+1)$-space $\CZ^{n+1}$, and their elementary
symmetric polynomials $
\sigma_k := \sum_{1 \leq i_1< \ldots <i_k \leq n+1}
\tilde{Z}_{i_1} \cdots \tilde{Z}_{i_k} $ for $
1 \leq k \leq n+1$. 
The ${\goth
A}_{n+1}$-invariant polynomials in
$\CZ[\tilde{Z}]$ are generated by the above
$\sigma_k$s and $
\delta := \prod_{i < j} (\tilde{Z}_i -
\tilde{Z}_j)$  
with a relation $
\delta^2= \widetilde{F}(\sigma_1, \sigma_2,
\ldots,
\sigma_{n+1})$
for certain polynomial $\widetilde{F}$. In fact,
$\widetilde{F}$  is a (quasi-)homogeneous
polynomial of degree $n(n+1)$ with the
weights of
$\sigma_k$ and 
$\delta$ equal to $k,
\frac{n(n+1)}{2}$ respectively. Denoted by 
$s_k, d$ the restriction
functions of $\sigma_k, \delta$ on $V$ respectively. Then
$s_1 $ is the zero function, and $V/{\goth
S}_{n+1} = \CZ^n$ via the coordinates $(s_2, \ldots,
s_{n+1})$. The orbifold $S_{{\goth
A}_{n+1}} (=V/{\goth A}_{n+1})$ is a double cover
of
$\CZ^n$, 
$$
S_{{\goth
A}_{n+1}} 
\longrightarrow \CZ^n = V/{\goth S}_{n+1}  \ .
$$
Then 
$V/{\goth S}_{n+1}$ can be realized as a
hypersurface in $\CZ^{n+1}$ with the
equation,
\be
S_{{\goth A}_{n+1}} : \ \ \
d^2 = F_n (  s_2,
\ldots, s_{n+1} ) \ , \ \ (d, s_2, \ldots, s_{n+1}
)
\in
\CZ^{n+1} \ , 
\ele(SA)
where $F_n (s_2, \ldots, s_{n+1}) :=
\widetilde{F}( 0, s_2,
\ldots, s_{n+1} )$. The 
polynomial $F_n(s_2, \ldots, s_{n+1}) $ has a
lengthy expression in general. Here we list the
polynomial
$F_n$ for $n=3, 4$:
\bea(lll)
F_3(s_2, s_3, s_4) & = & - 4s_2^3 s_3^2 -27 s_3^4
+ 16 s_2^4 s_4 -128 s_2^2s_4^2 + 144 s_2s_3^2 s_4
+ 256 s_4^3 
\ ; \\

F_4(s_2, s_3, s_4, s_5)& = & -4 s_2^3 s_3^2s_4^2 -
27 s_3^4 s_4^2 +16 s_2^4s_4^3 +144 s_2s_3^2s_4^3 
-128 s_2^2s_4^4  +256 s_4^5 
-72s_2^4s_3s_4s_5 

\\ 
& &
+108 s_3^5s_5
-630 s_2s_3^3s_4s_5 -1600s_3s_4^3s_5 +560
s_2^2s_3s_4^2s_5 +16 s_2^3 s_3^3s_5-900
s_2^3s_4s_5^2 \\
& &+ 2250s_3^2s_4s_5^2
+2000s_2s_4^2s_5^2 +108 s_2^5s_5^2 +825
s_2^2s_3^2s_5^2  -3750 s_2s_3s_5^3 + 3125 s_5^4 .
\elea(eqnS)

\section{${\bf A_1(4)}$-Singularity and Flop of
4-folds}
We now study the $A_1(n)$-singularity with $n\geq 4$. The set of $N$-integral elements in $\Delta$ are given by
$$
\bigtriangleup \cap N = \{ e^j \ | \ 1 \leq j
\leq n \}  \cup \{ v^{i,j} \ | \ 1
\leq i < j \leq n \} \ , 
$$
where $ v^{i,j} :=
\frac{1}{2}(e^i + e^j )$ for $i \neq j$.
Other than the simplex $\Delta$ itself, there is only one integral polytope decomposition of $\Delta$  invariant under all permutations of coordinates, and we will denote it by $\Xi$. $\Xi(n-1)$ consists of $n+1$ elements:
$\bigtriangleup_{i} \ ( 1 \leq i \leq n)$ and $\Diamond$, 
where $\bigtriangleup_{i}$ is the simplex generated by $e^i$ and
$v^{i,j}$ for $j \neq i$, and 
$\Diamond$ is the closure of 
$ \bigtriangleup\setminus \bigcup_{i=1}^n\bigtriangleup_{i}$, 
equivalently $\Diamond=$ the convex hull spanned by $v^{i, j}$s for
$i \neq j$. The lower dimensional polytopes of
$\Xi $ are the faces of those in $\Xi (n-1)$.  $X_{\Xi}$ has the trivial
canonical sheaf. For $n= 2, 3$, $X_{\Xi}$  is a crepant resolution of
$S_{A_1(n)}$.  For $n = 4$, one has the
following result.
\begin{lma}  \label{lem:SingXi}
For $n = 4$, the toric variety $X_{\Xi}$ is
smooth except one isolated singularity, which is
the 0-dimensional T-orbit corresponding to $\Diamond$.
\end{lma}
{\it Proof.} In general, for $n \geq 4$, it is easy to see that for each
$i$, the vertices of $\Delta_i$ form a $\ZZ$-basis of $N$, e.g., say $i=1$, it follows from $| A_1(n)| = 2^{n-1}$, and $
{\rm det}(e^1, v^{1,2}, \cdots, v^{1,n} ) =
\frac{1}{2^{n-1}}$. 
Hence $X_{\Xi}$ is non-singular near the $T$-orbits associated to simplices in $\Delta_i$. As  $\Diamond$ is not a simplex, ${\rm orb}(\Diamond)$ is always a singular point of $X_\Xi$. For $n=4$, the statement of smoothness of $X_{\Xi}$  except ${\rm orb}(\Diamond)$ follows from the fact that 
for $1 \leq i \leq 4$, the vertices $v^{i, j} (j \neq i)$ of 
$X_{\Xi}$, together with $\frac{1}{2} \sum_{j=1}^4 e^j$, form 
a $N$-basis.
$\Box$ \par \vspace{.2in} 

\begin{rmk}\label{rmk} For $n \geq 4$, the following properties hold for
0-dimensional $T$-orbits of $X_{\Xi}$.

(1) Denote $x_{\Delta_j} := {\rm
orb}(\Delta_j) \in X_{\Xi}$ for $1 \leq j \leq
n$. The inverse of the matrix spanned by vertices
of $\Delta_j$, $
( v^{1,j} , \cdots, v^{j-1,j}, e^j,  v^{j+1,
j}, \cdots, v^{n, j})^{-1} $, 
gives rise to affine coordinates $(U_1,
\ldots, U_n)$ centered at $x_{\Delta_j}$ such that  
$
U_i = Z_i^2 \ \ (i \neq j)$, and $ U_j= 
\frac{Z_j}{Z_1
\cdots \check{Z}_j \cdots Z_n}$. 
Hence $
I(x_{\Delta_j}) = \langle Z_j , Z_i^2 , i \neq j  \rangle + 
I (o)$, and we have
the regular $A_1(n)$-module structure of 
$\CZ[Z]/I(x_{\Delta_j})$,
\bea(l)
\CZ[Z]/I(x_{\Delta_j}) \simeq \bigoplus \{ \CZ Z^I \ | 
I=(i_1, \ldots, i_n) , \ i_j=0 , 
i_k =0, 
1 \ {\rm for } \ \ k
\neq j \} \ .
\elea(A1reg)

(2) We shall denote $x_{\Diamond}:= {\rm
orb}(\Diamond)$ in $X_{\Xi}$. The singular structure of
$x_{\Diamond}$ is determined by the $A_1(n)$-invariant polynomials corresponding to the $M$-integral elements in the cone dual to the
one generated by $\Diamond$ in $N_{\RZ}$. 
So the
$A_1(n)$-invariant polynomials are generated by  $
X_j := Z_j^2 $ and $ Y_j := \frac{Z_1 \cdots
\check{Z_j}
\cdots Z_n}{Z_j}$. 
Hence $
I(x_{\Diamond}) = \langle Z_1 \cdots
\check{Z_j}
\cdots Z_n  \rangle_{1 \leq j \leq n}  + 
I (o)$. 
Note that for $n=3$, the  $Y_j$s indeed form
the minimal generators for the invariant
polynomials, which implies the smoothness of
$X_{\Xi}$. For $n \geq 4$,  
$x_{\Diamond}$ is a singularity,  
not of the hypersurface type. For
$n=4$, the $X_j, Y_j \ (1 \leq j \leq 4)$ form a
minimal set of  generators of invariant
polynomials, hence the structure near
$x_{\Diamond}$ in $X_{\Xi}$ is the 4-dimensional
affine variety in $\CZ^8$  defined by the
relations:
\begin{eqnarray}
x_iy_i= x_j y_j ,
\ \
\  x_ix_j = y_{i^\prime}y_{j^\prime} \ , \ (x_i, y_i)_{1
\leq i \leq 4} \in \CZ^8 , 
\label{4sing}
\end{eqnarray}
where $i \neq j$ with $\{ i^\prime, j^\prime \} $
the complementary pair of $\{i, j\}$.
\end{rmk}

\begin{figure}[ht]
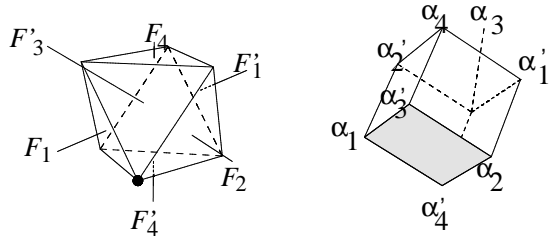

\hfil{\psfig{figure=./octeface.ps}}\hskip 0.5in {\psfig{figure=./cubic.ps}}\hfil

\caption {Dual pair of octahedron and cube: Faces
$F_j, F_j^\prime$ of octahedron  dual
to vertices $\alpha_j, \alpha_j^\prime$ of cube. The face of the cube in gray color
corresponds to the dot ``$\bullet$" in the
octahedron.}
\label{fig:octe}
\end{figure}

For the rest of this section, we shall consider only the case $n=4$. We are going to discuss the structure of $\hl^{A_1(4)}(\CZ^4)$ and its
connection with crepant resolutions of $S_{A_1(4)}$. The simplex 
$\Delta$ is a tetrahedron, and $\Diamond$ is an octahedron; both are acted on by the symmetric group ${\goth S}_4$. The dual polygon of $\Diamond$ is the cube. The facets of the octahedron $\Diamond$ are labeled by 
$F_j, F_j^\prime$ for $1 \leq j \leq 4$, where $
F_j = \Diamond \cap \bigtriangleup_j$ and $
F_j^\prime = 
\{ \sum_{i=1}^4 x_ie^i \in \Diamond \ | \ x_j =0
\}$. 
The dual of $F_j$, $F_j^\prime$ in the cube are vertex, denoted by $\alpha_j$, $\alpha_j^\prime$ as in Fig. \ref{fig:octe}.

\begin{figure}[ht]
\begin{center}
\mbox{\psfig{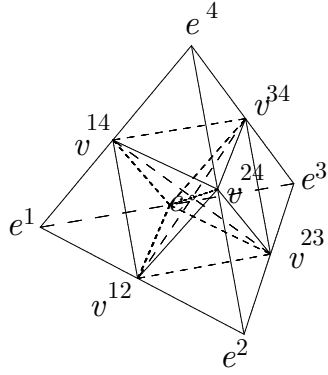}}
\caption{The rational simplicial decomposition
$\Xi^*$ of $\Delta$ for $n=4, r=1$.}  
\label{fig:tatrhilb}
\end{center}
\end{figure}

Consider the rational simplicial decomposition $\Xi^*$ of
$\Delta$, which is a refinement of $\Xi$ by adding the center 
$c := \frac{1}{4} \sum_{j=1}^4 e^j $
as a vertex with the barycentric decomposition
of $\Diamond$ in $\Xi$, (see Fig. \ref{fig:tatrhilb}). Note that
$c \not\in N$ and $2c \in N$.
For convenience, we shall use the following convention: 
\par \noindent 
{\bf Notation.} Let $G$ be a diagonal group acting
on $\CZ[Z]$. Two monomials $m_1, m_2$ in
$\CZ[Z]$  are said to be
$G$-equivalent, denoted by $m_1{\buildrel G \over
\sim } m_2$ or simply by $m_1\sim m_2$ , if 
$m_1/m_2$ is a $G$-invariant function. 

\begin{thm}\label{th:A1(4)}
For $G= A_1(4)$, we have $
\hl^G(\CZ^4)  \simeq 
X_{\Xi^*}$, 
which is non-singular with the canonical bundle 
$\omega = {\cal O}_{X_{\Xi^*}} (E)$, where $E$ is
an irreducible divisor isomorphic to the triple
product of
$\PZ^1$, 
\begin{eqnarray}
E = \PZ^1 \times \PZ^1 \times \PZ^1 \ . 
\label{eq:E}
\end{eqnarray}
Furthermore for $\{ i, j, k \} = \{1,2,3
\}$, the normal bundle of $E$ when restricted on the fiber $\PZ^1_k \ (\simeq \PZ^1)$, for the projection $E$ to $\PZ^1\times \PZ^1$ via the
$(i, j)$-th factor,
\begin{eqnarray}
p_k : E \longrightarrow \PZ^1 \times \PZ^1 \
, 
\label{eq:proj}
\end{eqnarray}
is the
$(-1)$-hyperplane bundle:
\begin{eqnarray}
{\cal O}_{X_{\Xi^*}}(E) \otimes {\cal
O}_{\PZ^1_k}
\simeq  {\cal O}_{\PZ^1}(-1)  \ . 
\label{eq:f(-1)}
\end{eqnarray}
\end{thm}
{\it Proof.}  First we show the smoothness of the toric
variety $X_{\Xi^*}$. The octahedron $\Diamond$ of $\Xi$ is decomposed into eight simplices of $\Xi^*$ corresponding to faces $F_j, F_j^\prime$ of $\Diamond$. Denote $C_j$ (resp. $C_j^\prime$) the simplex of $\Xi^*$ spanned by
$c$ and $F_j$ (resp. $F_j^\prime$); $x_{C_j}, x_{C_j^\prime}$ are the corresponding 0-dimensional $T$-orbits in $X_{\Xi^*}$. The smoothness of affine space in $X_{\Xi^*}$ near $x_{C_j}, x_{C_j^\prime}$ follows from the $N$-integral criterion of the cones in $N_{\RZ}$ generated by $C_j, C_j^\prime$. The coordinate system is given by the integral basis of $M$
which generates the cone dual to the cone spanned by $C_j$ ( $C_j^\prime$). As examples, for $C_1, C_2^\prime$, the coordinates are determined by the row vectors of the following square matrix:

$$
\begin{array}{ll}
{\rm cone}(C_1)^* , & {\rm cone}(C_2^\prime)^*  \\
(2c, v^{1,2}, v^{1,3},
v^{1,4})^{-1} = 
\left( \begin{array}{cccc}
-1 &   1 & 1 & 1  \\
1 & 1 & -1 & -1 \\
1 & -1 & 1 & -1 \\
1 & -1 & -1 & 1
\end{array} \right) \ , & 
(v^{3,4}, 2c, v^{1,4}, v^{1,3})^{-1} = 
\left( 
\begin{array}{cccc}
-1 & -1 & 1 & 1  \\
0 & 2 & 0 & 0 \\
1 & -1 & -1 & 1 \\
1 & -1 & 1 & -1
\end{array} \right) \ .
\end{array}
$$
The coordinate functions of
$X_{\Xi^*}$ centered at $x_{C_1}$ are given by $
(U_1, U_2, U_3, U_4)= (\frac{Z_2Z_3Z_4}{Z_1}, \frac{Z_1Z_2}{Z_3Z_4}, \frac{Z_1Z_3}{Z_2Z_4}, \frac{Z_1Z_4}{Z_2Z_3})$ with $
I (x_{C_1}) = \langle Z_2Z_3Z_4, Z_1Z_2, Z_1Z_3, 
Z_1Z_4\rangle + I (o)$,
and the coordinates near $x_{C_2^\prime}$ are
$
(U_1^\prime, U_2^\prime, U_3^\prime, U_4^\prime)=(  \frac{Z_3Z_4}{Z_1Z_2}, Z_2^2 ,  \frac{Z_1Z_4}{Z_2Z_3},  \frac{Z_1Z_3}{Z_2Z_4}) $ 
with $
I (x_{C_2^\prime}) = \langle Z_3Z_4,Z_2^2, Z_1Z_4, 
Z_1Z_3 \rangle + I(o)$. 
By the Remark \ref{rmk} (1), one has the smooth coordinate system centered at $x_{\Delta_j}$ in 
$X_{\Xi^*}$. For $\Delta_1$, by
$$
\begin{array}{ll}
{\rm cone}(\Delta_1)^* , & (e^1, v^{1,2}, v^{1,3},
v^{1,4})^{-1} = 
\left( \begin{array}{cccc}
1 &  -1 & -1 & -1  \\
0 & 2 & 0 & 0 \\
0 & 0 & 2 & 0 \\
0 & 0 & 0 & 2
\end{array} \right) \ , 
\end{array}
$$
one has the coordinate system near
$x_{\Delta_1}$,
$
(V_1, V_2, V_3, V_4)  =(\frac{Z_1}{Z_2 Z_3 Z_4} , Z_2^2,
 Z_3^2 , Z_4^2 )$ 
with $
I (x_{\Delta_1}) = \langle Z_1, Z_2^2, Z_3^2, 
Z_4^2 \rangle + I (o)$. 
Now we are going to show that $\CZ[Z]/I(y)$ is a regular $G$-module for $y \in  X_{\Xi^*}$. For an element $y$ in the affine neighborhood of
$x_{\Delta_1}$ with the coordinates 
$V_i=v_i,( 1 \leq i\leq 4) $, one has
\begin{equation} \label{eq:IdDelta}
I(y) = \langle Z_1-v_1 Z_2Z_3Z_4,
Z_2^2-v_2 ,  Z_3^2-v_3 ,  Z_4^2-v_4  \rangle 
\end{equation}
The set of monomials, $
\left\{1,Z_2,Z_3,Z_4,Z_2Z_3,Z_2Z_4,Z_3Z_4,
Z_2Z_3Z_4\right\}$, 
gives rise to a basis of $\CZ[Z]/I(y)$ for $v_i\in\CZ$; hence $\CZ[Z]/I(y)$ is a regular $G$-module. For $y$ near $x_{C_1}$ with the coordinates $U_i=u_i, ( 1 \leq i \leq 4)$, we have 
\begin{equation}\label{eq:IdC}
I (y) = \langle Z_2Z_3Z_4-u_1Z_1, 
Z_1Z_2-u_2Z_3Z_4, Z_1Z_3-u_3Z_2Z_4,
Z_1Z_4-u_4Z_2Z_3\rangle + I_{G}(y) \ ,
\end{equation}
where $I_{G}(y)=\langle Z_1Z_2Z_3Z_4-u_1^2u_2u_3u_4,
Z_1^2-u_1u_2u_3u_4, Z_2^2-u_2u_1, Z_3^2-u_3u_1,
Z_4^2-u_4u_1 \rangle $. 
This implies that $\CZ[Z]/I(y)$ is a regular
$G$-module with a basis represented by $
\left\{ 1,Z_1,Z_2,Z_3,Z_4,
Z_2Z_3,Z_3Z_4,Z_2Z_4 \right\}$.  
Similarly, the same conclusion holds for $y$ near 
$x_{C_2^\prime}$ with the coordinates
$U_i^\prime=u_i^\prime, ( 1 \leq i \leq 4)$, in
which case we have 
\begin{equation}\label{eq:IdCp}
I (y) = \langle Z_3Z_4-u_1^\prime Z_1Z_2, Z_1Z_4-u_3^\prime Z_2Z_3, Z_1Z_3-u_4^\prime Z_2Z_4\rangle + I_{G}(y) \ ,
\end{equation}
with $I_{G}(y)=\langle Z_1Z_2Z_3Z_4-
{u_2^\prime }^2u_1^\prime u_3^\prime u_4^\prime ,
Z_2^2-u_2^\prime , Z_1^2-u_2^\prime
u_3^\prime u_4^\prime , Z_3^2-u_1^\prime
u_2^\prime u_4^\prime , Z_4^2-u_1^\prime
u_2^\prime u_3^\prime  \rangle$, and a basis of
$\CZ[Z]/I(y)$ represented by $
\left\{1,Z_1,Z_2,Z_3,Z_4,Z_1Z_2,Z_2Z_3,Z_2Z_4
\right\}$. 
The same argument can equally be applied to all 
affine charts centered at $x_{\Delta_j},
x_{C_j}, x_{C_j^\prime}$. Therefore we obtain a morphism
$$
\lambda : X_{\Xi^*} \longrightarrow {\rm
Hilb}^G(\CZ^4) \ , \ \ {\rm with} \ \ 
I(\lambda(y))=I(y),  \ y\in X_{\Xi^*} \ . 
$$

We are going to show that the above morphism
$\lambda$ is an isomorphism by
constructing its inverse morphism. Let $y'$ be an element of $ \hl^{G}(\CZ^4)$, 
represented by  a $G$-invariant ideal $J\subset
\CZ[Z]$ with $\CZ[Z]/J $ as the regular $G$-module.  By Gr\"{o}bner basis techniques \cite{CLO}, for a given  monomial order,  there is 
a monomial ideal ${\rm lt}(J)$, 
consisting of all leading monomials of elements in $J$, such
that  the monomial base of $\CZ[Z]/{\rm lt}(J)$
also gives rise to a basis of $\CZ[Z]/J$. By this fact, we
shall first determine the $G$-invariant monomial ideal $J_0$ in ${\rm
Hilb}^G(\CZ^4)$. For a monomial $I$, we shall denote $I^\dagger$  
the set of monic monomials not in $I$.  Since all 
nonconstant $G$-invariant monomials are in
$J_0$, we have ${Z_j}^2, Z_1Z_2Z_3Z_4 \in J_0$. 
Hence  $\eM$ is contained in the set 
$\bas : = \{ Z^I \ | \ I=(i_1,..,i_4), \
i_1i_2i_3i_4=0, i_j\le 1 \} $. For a nontrivial character
$\rho$ of $G$, the $\rho$-eigenspace of
$I(o)^\bot$ for the element $o \in S_G$ is of 
dimension 2. This implies that
for $m_1
\in \bas$ not equal to $ 1$, there exists exactly one
$m_2\in\bas$ not equal to $ m_1$ with 
$m_2 \sim m_1$. When $J_0=I(x_{\Delta_1})$,
$I(x_{\Delta_1})^\bot$ has a monomial basis $W:=I(x_{\Delta_1})^\dagger$ consisting of eight elements ${Z_2}^{i_2}{Z_3}^{i_3}{Z_4}^{i_4}$, $0\le i_j\le 1$, 
and they form a basis of the $G$-regular representation. 
By replacing some monomials in $W$ by the other $G$-equivalent ones in $\bas$, one obtains a $G$-regular basis $W^\prime$. Denote $W_0$ the set of monic monomials in $\CZ[Z]$. The $W^\prime$s satisfying $W_0 \cdot (W_0-{W^\prime}) \subset(W_0-{W^\prime})$ are in one-to-one correspondence with monomial ideals $J_0$s in $ {\rm
Hilb}^{G}(\CZ^4)$ by  the relation $J_0=\langle W_0-{W^\prime}\rangle_\CZ$, hence $W^\prime  = J_0^\dagger$. By direct counting, there are 
twelve such $W^\prime$ and the corresponding twelve $J_0$'s,
are exactly those $I(x_\RT)$ for $\RT \in \Xi^*(3)$. The correspondence between
$W^\prime$ and $\RT$ by the relation $W^\prime = 
I(x_\RT)^\dagger $ is given as follows:
\bea(cll)
\left\{1,Z_2,Z_3,Z_4,Z_2Z_3,Z_2Z_4,Z_3Z_4,
Z_2Z_3Z_4\right\} & \leftrightarrow &{\Delta_1},
\\
\left\{1, Z_1, Z_3,Z_4,Z_1Z_4,Z_1Z_3,Z_3Z_4,
Z_1Z_3Z_4 \right\} &
\leftrightarrow & {\Delta_2}, \\
\left\{1,Z_1, Z_2,Z_4,Z_1Z_4,Z_2Z_4,Z_1Z_2,
Z_1Z_2Z_4 \right\} &
\leftrightarrow &{\Delta_3}, \\
\left\{1,Z_1,Z_2,Z_3,Z_2Z_3,Z_1Z_3,Z_1Z_2,
Z_1Z_2Z_3 \right\} &
\leftrightarrow &{\Delta_4}, 
\\
\left\{1,Z_1,Z_2,Z_3,Z_4,Z_2Z_3,Z_2Z_4,Z_3Z_4
\right\} &
\leftrightarrow &{C_1}, \\
\left\{1,Z_1 ,Z_2,Z_3,Z_4,Z_1Z_4,Z_1Z_3,Z_3Z_4
\right\} & 
\leftrightarrow &{C_2}, \\
\left\{1,Z_1,Z_2,Z_3,Z_4,Z_1Z_4,Z_2Z_4,Z_1Z_2
\right\} &
\leftrightarrow & {C_3}, \\
\left\{1,Z_1,Z_2,Z_3,Z_4,Z_2Z_3,Z_1Z_3,Z_1Z_2
\right\} &
\leftrightarrow &{C_4}, \\
\left\{1,Z_1,Z_2,Z_3,Z_4,Z_1Z_4,Z_1Z_3,Z_1Z_2
\right\} &
\leftrightarrow & {C_1^\prime}, \\
\left\{1,Z_1,Z_2,Z_3,Z_4,Z_2Z_3,Z_2Z_4,Z_1Z_2
\right\} &
\leftrightarrow & {C_2^\prime}, \\
\left\{1,Z_1,Z_2,Z_3,Z_4,Z_2Z_3,Z_1Z_3,Z_3Z_4
\right\} &
\leftrightarrow &{C_3^\prime}, \\
\left\{1,Z_1,Z_2,Z_3,Z_4,Z_1Z_4,Z_2Z_4,Z_3Z_4
\right\} &
\leftrightarrow & {C_4^\prime}.
\elea(A2)
Now we consider an ideal $J$ in $\CZ[Z]$ which defines an element of 
$\hl^G(\CZ^4)$. By the Gr\"{o}bner basis 
argument as before, there is a monomial ideal $J_0(={\rm lt}(J))$ 
such that $\eM$ gives rise to a basis of $\CZ[Z]/J$, and  $J_0=I(x_\RT)$ for some $\RT\in \Xi^*(3)$.  
For $p\in\CZ[Z]$, the element $p+J\in\CZ[Z]/J$ is uniquely expressed in the
form, $p+J = \sum_{m\in\eM}\gamma (p)_m m+J$,
i.e., 
$p-\sum_{m\in\eM}\gamma (p)_m m\in J$. In particular, for a monomial  $p$ in $\CZ[Z]$, ( i.e., $p \in W_0$), we have $g\cdot (p-\sum_{m\in\eM}\gamma (p)_m m) \in J$ for $g\in G$. This implies $p-\sum_{m\in\eM}\gamma
(p)_m \mu_g(p)^{-1}\mu_g(m) m \in J$,  where 
$\mu_g(m),\mu_g(p)\in \CZ^*$ are the 
the character values of $g$ on $m, p$
respectively;  hence 
$$
\sum_{m\in\eM}\gamma (p)_m
\left[{\mu_g(p)}^{-1}\mu_g(m)-1\right]m \in J .
$$
As $\eM$ represents a $G$-regular basis for
$\CZ[Z]/J$, we have
$\gamma(p)_m\left[{\mu_g(p)}^{-1}\mu_g(m)-1\right]=0$ 
for $p \in W_0$, $m\in \eM$ and  $g\in G$. Furthermore, for each $p \in W_0$, there exists an unique element, denoted by $p_{\eM}$, in $\eM$ with the property $p \sim p_{\eM}$. Hence for $m\in \eM$, $ m\neq
p_{\eM}$ if and only if   
$\left[{\mu_g(p)}^{-1}\mu_g(m)-1\right]\ne 0$ for
some $g \in G$, in which case $\gamma(p)_m = 0$.
Therefore $p - \gamma (p)_{p_{\eM}}
p_{\eM} \in J$, and $J$ is the ideal with the
generators: 
\bea(l)
J = \langle \ p -\gamma
(p)_{p_{\eM}} p_{\eM} \ | \ p \in W_0 \cap J_0 \
\rangle
\ .
\elea(Jgen)
Indeed in the above expression of $J$, it suffices to consider those $p$s which from a minimal set of monomial generators of $J_0$.   
Now we are going to assign an element of $X_{\Xi^*}$ for a given $J \in \hl^G(\CZ^4)$. If the monomial ideal $J_0$ associated to 
$J$ in our previous discussion is equal to $I(x_{C_1})$, a minimal set of monomial generators of $J_0$ and the basis representative set $J_0^{\dagger}$ of $\CZ[Z]/J$ are given by
$$J_0=\langle Z_1^2,Z_2^2,..,Z_4^2,Z_1Z_2,Z_1Z_3,Z_1Z_4,Z_2Z_3Z_4\rangle,$$
$$J_0^{\dagger}=\{1,Z_1,Z_2,Z_3,Z_4,Z_2Z_3,Z_2Z_4,Z_3Z_4\}. $$
By (\ref{Jgen}),  $J$ contains the ideal generated by
$p - \gamma(p)_{p_{\eM}} p_{\eM}$ for
$p=Z_i^2, Z_1Z_2, Z_1Z_3, Z_1Z_4, Z_2Z_3Z_4$ for $1 \leq i \leq 4$, 
which has the
colength at most $8$ in $\CZ[Z]$. Therefore 
one obtains
$$\begin{array}{lll}
J =&\langle 
Z_1Z_4-\gamma_{14}Z_2Z_3,\ 
Z_1Z_3-\gamma_{13}Z_2Z_4,\
Z_1Z_2-\gamma_{12}Z_3Z_4, Z_2Z_3Z_4-\gamma_{234}Z_1,  \\  
{}&{Z_1}^2-\gamma_{1},\ 
{Z_2}^2-\gamma_{2},\ {Z_3}^2-\gamma_{3},\
{Z_4}^2-\gamma_{4}\ 
\rangle 
\end{array}
$$
Moreover, by 
$$
\begin{array}{l}
0 \equiv Z_2({Z_1}^2-\gamma_{1})-
Z_1(Z_1Z_2-\gamma_{12}Z_3Z_4)
\equiv(\gamma_{12}\gamma_{13}\gamma_4-\gamma_1)Z_2
\ {\rm  (mod } J),
\end{array}
$$  
and $Z_2 \in {J_0}^{\dagger}$, one has
$$
\gamma_1=\gamma_{12}\gamma_{13}\gamma_4 \ \ \ .
$$
By
$$
\begin{array}{l}
0 \equiv Z_1({Z_4}^2-\gamma_{4})-
Z_4(Z_1Z_4-\gamma_{14}Z_2Z_3)
\equiv(\gamma_{14}\gamma_{234}-\gamma_4)Z_1
\ {\rm  (mod } J),
\end{array}
$$  
one obtains
$$
\gamma_2=\gamma_{234}\gamma_{12}.
$$
Similarly, one has
$$
\gamma_3=\gamma_{234}\gamma_{13},
\gamma_4=\gamma_{234}\gamma_{14}, \ \ \ 
$$
Therefore, all $\gamma_I$s are expressed as functions of $\gamma_{12}, \gamma_{13}, \gamma_{14},
\gamma_{234}$. This implies $J=I(y)$ for an element
$y$ of $X_{\Xi^*}$ in the affine neighborhood
$x_{C_1}$ with the coordinate $(U_i=u_i)$ 
by the relations, 
$$
u_1=\gamma_{234}, \ \ \
u_2=\gamma_{12}, \ \ \
u_3=\gamma_{13}, \ \ \ 
u_4=\gamma_{14}.
$$
The above $y$ is defined
to be the element $\lambda^{-1}(J)$ in $X_{\Xi^*}$
for the ideal
$J$ under the inverse map of
$\lambda$. The method can equally be applied to
ideals  $J$ associated
to another monomial ideal $J_0$.  

For $J_0= I(x_{C_2^\prime})$, we have
$$
J= \langle Z_1Z_3-\gamma^\prime_{13}Z_2Z_4, \ 
Z_1Z_4-\gamma^\prime_{14}Z_2Z_3, \
Z_3Z_4-\gamma^\prime_{34}Z_1Z_2,\ {Z_1}^2-\gamma^\prime_{1}, \ 
{Z_2}^2-\gamma^\prime_{2}, \
{Z_3}^2-\gamma^\prime_{3}, \ {Z_4}^2-\gamma^\prime_{4}\rangle.
$$
We claim that the variables 
$\gamma_2^\prime$, $\gamma_{34}^\prime$, $\gamma_{13}^\prime$,  
$\gamma_{14}^\prime$ form a system of coordinates near $I(x_{C_2^\prime})$, i.e., all the $\gamma_I^\prime$s can be 
expressed as certain polynomials of these four values. Indeed, we are going to show  
$\gamma^\prime_1= \gamma^\prime_2\gamma^\prime_{13}\gamma^\prime_{14},$
$\gamma^\prime_3= \gamma^\prime_2\gamma^\prime_{13}\gamma^\prime_{34}$ 
and 
$\gamma^\prime_4= \gamma^\prime_2\gamma^\prime_{14}\gamma^\prime_{34}.$\footnote{Note that the group $G$ in Section 6.1 of \cite{N} (page 777) is the $A_1(4)$ of Theorem \ref{th:A1(4)} in this paper. However, we would consider that the statement in \cite{N} about the singular property of $\hl^G (\CZ^4)$  by using the structure of $I(\Gamma_3)(u)$ there, is not correct. Indeed, by 
identifying $Z_2, Z_3, Z_4, Z_1$ with $x, y, z, w$, and 
$\gamma_2^\prime,  \gamma_3^\prime,  \gamma_4^\prime$,  
$\gamma_1^\prime, \gamma_{34}^\prime, \gamma_{13}^\prime$, 
$\gamma_{14}^\prime$ with $u_1, u_2, \cdots , u_7$ respectively, the ideal $J$ in our discussion corresponds to $I(\Gamma_3)(u)$ in \cite{N}. Then through the three relations we have obtained here,  one can easily verify that all the relations among the $u_i$s listed in page 778 of \cite{N} hold.}
By 
$$Z_1(Z_1Z_4-\gamma^\prime_{14}Z_2Z_3)- Z_4({Z_1}^2-\gamma^\prime_1) = 
-\gamma^\prime_{14}Z_1Z_2Z_3+\gamma^\prime_1Z_4\in J,
$$
one has
$$Z_2(-\gamma^\prime_{14}Z_1Z_2Z_3+\gamma^\prime_1Z_4)+\gamma^\prime_{14}
Z_1Z_3(Z_2^2-\gamma_2^\prime) = \gamma^\prime_1Z_2Z_4- 
\gamma^\prime_{14}\gamma_2^\prime Z_1Z_3\in J , $$ 
hence 
$$(\gamma^\prime_1Z_2Z_4- \gamma^\prime_{14}\gamma_2^\prime 
Z_1Z_3)+ \gamma^\prime_{14}\gamma_2^\prime (Z_1Z_3-\gamma_{13}^\prime 
Z_2Z_4)=(\gamma^\prime_1-\gamma_2^\prime \gamma_{13}^\prime 
\gamma^\prime_{14})Z_2Z_4 \in J .
$$
By the description in (\req(A2)) for $C_2^\prime$, $Z_2Z_4$ is an element in $\eM$, hence 
represents a basis element of 
$\CZ[Z]/J$. The relation $(\gamma^\prime_1-\gamma_2^\prime \gamma_{13}^\prime 
\gamma^\prime_{14})Z_2Z_4 \in J$ implies 
$$\gamma^\prime_1-\gamma_2^\prime \gamma_{13}^\prime 
\gamma^\prime_{14}= 0 \ .
$$ 
By interchanging the indices $1$ and $3$, (resp. $1$ and $4$), in the above 
derivation and regarding $\gamma^\prime_{i j} = \gamma^\prime_{j i}$, we obtain
$\gamma^\prime_3= \gamma^\prime_2\gamma^\prime_{13} \gamma^\prime_{34}$ 
(resp. 
$\gamma^\prime_4= \gamma^\prime_2\gamma^\prime_{14}\gamma^\prime_{34}$).
 Thus, $\gamma^\prime_2$,  $\gamma^\prime_{13}$, $\gamma^\prime_{13}$ and $\gamma^\prime_{34}$ form the four independent parameters to describe  
the ideals $J$ near $J_0= I(x_{C_2^\prime})$ with the regular $G$-module 
$\CZ[Z]/J$. Therefore  
$J=I(y)$ for  $y$ near
$x_{C_2^\prime}$ with the coordinates
$(U^\prime_i=u^\prime_i)$ via the relations,
$$
u_2^\prime=\gamma^\prime_2, \ \ 
u_1^\prime=\gamma^\prime_{34}, \ \
u_3^\prime=\gamma^\prime_{14}, \ \ 
u_4^\prime=\gamma^\prime_{13}.
$$

For $J_0=I(x_{\Delta_1})$, we have $
J= \langle Z_1-\gamma^{\prime\prime}_1 Z_2Z_3Z_4 ,
{Z_2}^2-\gamma^{\prime\prime}_2,
{Z_3}^2-\gamma^{\prime\prime}_3,
{Z_4}^2-\gamma^{\prime\prime}_4 \rangle$. 
Hence $J=I(y)$ for $y$ near $x_{\Delta_1}$ with the coordinates $(V_i=v_i)$ and the relations,
$v_i=\gamma^{\prime\prime}_i$ for
$1 \leq i \leq 4$. The previous discussions of
three cases can be applied to each of the twelve
monomial ideals $J_0$'s by a suitable change 
of indices. Hence one obtains an element $\lambda^{-1}(J)$ in $X_{\Xi^*}$ of an
ideal $J \in \hl^G(\CZ^4)$. 

However, one
has to verify  the correspondence  
$\lambda^{-1}$ so defined to be a single-valued map, namely,
for a given $J$ with two possible choices of $J_0$, the elements in
$X_{\Xi^*}$ assigned to $J$ through the previous
procedure through these two
$J_0$ are the same one. For example,
say $J= I(y_1)=I(y_2)$ for
$y_1$ near $x_{\Delta_1}$ with $(V_i=v_i)$, and
$y_2$ near $x_{C_1}$ with $ (U_i=u_i)$.
By (\ref{eq:IdDelta}), (\ref{eq:IdC}), both $
Z_2Z_3Z_4-u_1Z_1$ and $ Z_1-v_1Z_2Z_3Z_4$ are elements in $J$.
We claim that $u_1\ne 0$. Otherwise, both $Z_1$ and 
$Z_2Z_3Z_4$ are elements in $J$  with the
same $G$-character $\kappa$. Then the
$\kappa$-eigenspace in
$\CZ[Z]/J$ is the zero space, a contradiction to
the regular $G$-module property of $\CZ[Z]/J$. 
Hence one has $Z_1-{u_1}^{-1}Z_2Z_3Z_4\in J$, hence
$(v_1-{u_1}^{-1})Z_2Z_3Z_4\in J$.  
As $J=I(y_1)$ with $y_1$ near $x_{\Delta_1}$, $Z_2Z_3Z_4$ represents a  basis element of
$\CZ[Z]/J$. Hence $v_1={u_1}^{-1}$. 
By 
$Z_1Z_2-u_2Z_3Z_4 ,\ {Z_2}^2-v_2 \in J$, one has 
$v_2Z_1-u_2Z_2Z_3Z_4 (=
(Z_1Z_2-u_2Z_3Z_4)Z_2-({Z_2}^2-v_2)Z_1) \in J$. As
$Z_2Z_3Z_4 \not\in J$, one has  $u_2=0$ if 
$v_2=0$. 
When $v_2\ne 0$, we have, $
Z_1-u_2{v_2}^{-1}Z_2Z_3Z_4 \in J$, hence 
$$
(v_1-u_2{v_2}^{-1})Z_2Z_3Z_4\in J , \ \
\ u_2=v_1v_2
\ .
$$
Using the same argument, one can derives $u_j=v_1v_j$ for
$j=2,3,4$. These three relations, together with $u_1={v_1}^{-1}$, imply  
$y_1 = y_2$ in $X_{\Xi^*}$.

For $y_2$ near
$x_{C_1}$ with $(U_i=u_i)$, and 
$y_3$ near 
$x_{C_2^\prime}$ with $(U_i^\prime=u_i^\prime)$,
by (\ref{eq:IdC})  (\ref{eq:IdCp}), both 
$Z_1Z_2-u_2Z_3Z_4$ and $Z_3Z_4-u_1^\prime Z_1Z_2$ are elements in
$J$; furthermore, $u_2 , u_1^\prime $ are
non-zero by the fact that only one of $Z_1Z_2,
Z_3Z_4$ could be an element of $J$. By an argument similar
to the one before, one can show 
$$
u_1^\prime={u_2}^{-1}, \ \
u_3=u_4^{\prime}, \ \ u_4=u_3^{\prime}.
$$
By
$Z_2Z_3Z_4-u_1Z_1,\ Z_3Z_4-{u_1^\prime}Z_1Z_2,\
{Z_2}^2-{u_2^\prime}\in J$, we have
$$
(Z_2Z_3Z_4-u_1Z_1)Z_2\equiv (u_1^\prime
u_2^\prime-u_1)Z_1Z_2\equiv 0 \ \  {\rm mod }J \
. 
$$
As 
$Z_1Z_2$ represents a basis element of $\CZ[Z]/J$, one has $u_1= u_1^\prime u_2^\prime$. 
The four relations between $u_i$s and $u^\prime_i$s imply 
$y_2=y_3$ in $X_{\Xi^*}$. In this way, one can show directly that for a
given ideal
$J$ with
$J= I(y)=I(y')$ for $y, y'$ 
in $X_{\Xi^*}$, the elements $y$ and $y'$ are the same one by the relations of toric coordinates
centered at two distinct
$x_{\RT}$s. Hence we have obtained a well-defined
morphism $\lambda^{-1}$ from $\hl^G(\CZ^4)$ to
$X_{\Xi^*}$, then   
$\hl^G(\CZ^4)\simeq X_{\Xi^*}$. By  (\ref{K}), the 
canonical bundle of
$X_{\Xi^*}$.  is given by $
\omega  = {\cal O}_{X_{\Xi^*}}(E)$, 
where $E$ denotes the toric divisor $D_c$, which
is a 3-dimensional complete toric variety with
the toric data described by the star of $c$ in
$\Xi^*$, which is represented by
the octahedron in Fig. \ref{fig:octe}, where the cube in
Fig. \ref{fig:octe} represents the toric orbits' structure. 
Therefore
$E$ is isomorphic to the triple product of
$\PZ^1$ as in (\ref{eq:E}). The description of the
normal bundle of $E$ restricting on each
$\PZ^1$-fiber will follow by the direct
computation in toric geometry. For
example, for the fibers over the projection of $E$
onto $(\PZ^1)^2$ corresponding to the
2-convex set spanned by $v^{1,2}, v^{1,3},
v^{3,4}$ and $ v^{2,4}$, one can perform the
computation as follows. Let $(U_1, U_2, U_3,
U_4)$ be the local coordinates near
$x_{C_4^\prime}$ dual to the
$N$-basis 
$(2c, v^{1,2}, v^{1,3},  v^{2,3})$, and
let $(W_1, W_2, W_3,W_4)$ be the local coordinates 
near $x_{C_1}$ dual to $(2c, v^{1,2},
v^{1,3}, v^{1,4})$. By $2c= v^{1,4}+v^{2,3}$,
one has the relations, $
U_1=W_1W_4$ , $U_4= W_4^{-1}$, $ 
U_2=W_2$, $ U_3=W_3$. 
This shows that the restriction of the normal
bundle of $E$ on each fiber $\PZ^1$ over $(U_2,
U_3)$-plane is the $(-1)$-hyperplane bundle.
$\Box$ \par \vspace{.2in} \noindent
Note that the vector bundle ${\cal F}_{X_{\Xi^*}}$ over $X_{\Xi^*}$ 
in Theorem
\ref{th:A1(4)} carries the regular
$G$-module structure on each fiber with the local
frame of the vector bundle  provided by
the structure of $\CZ[Z]/I(x_{\RT})$ for $\RT
\in \Xi^*(3)$ with the representative in the list 
(\req(A2)). 
 
\begin{figure}[ht]
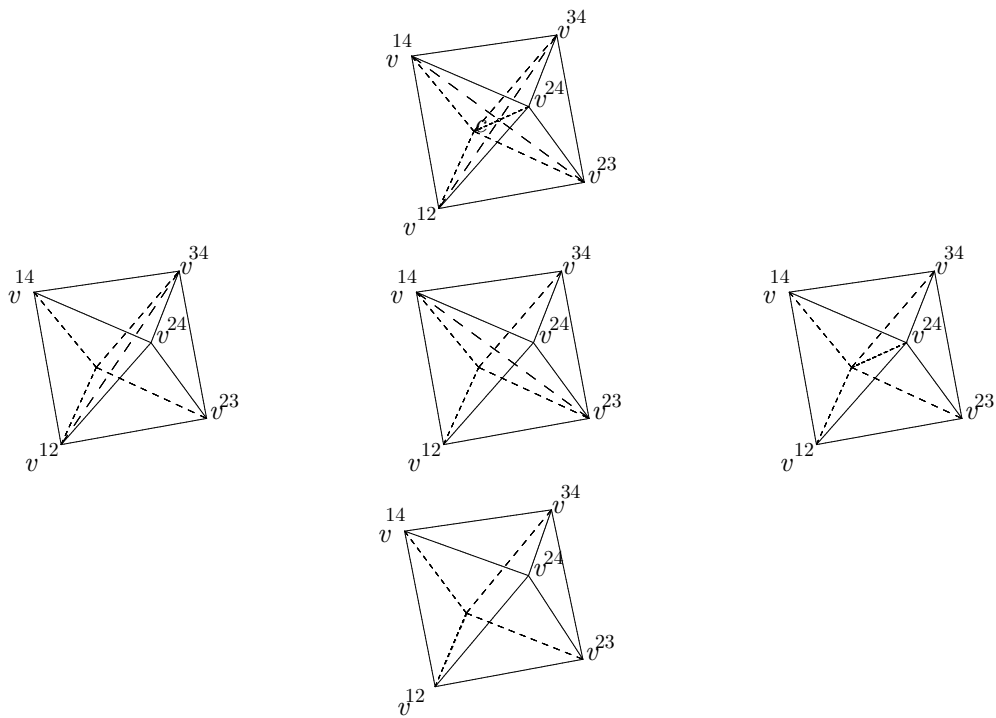

\begin{center}
\mbox{\psfig{figure=hilbs.ps}}

\mbox{\psfig{figure=crep1.ps}} {\hskip 50pt} 
\mbox{\psfig{figure=crep2.ps}}{\hskip 50pt} 
\mbox{\psfig{figure=crep3.ps}}

\mbox{\psfig{figure=octa.ps}}
\label{fig:hco}
\caption{Toric representation of 4-dimensional 
flops in the second row over a common singular base in the third row and dominated by the same 4-fold in the first row.}    
\end{center}
\end{figure}

By the standard
blowing-down criterion of an exceptional divisor,
the property (\ref{eq:f(-1)}) ensures the
existence of a smooth $4$-fold $(X_{\Xi^*})_k$ by
blowing-down the $\PZ^1$-family  along the
projection $p_k$ (\ref{eq:proj}) for each $k$. In
fact,
$(X_{\Xi^*})_k$ is also a toric variety
$X_{\Xi_k}$ with $\Xi_k$ defined by the refinement
of $\Xi$ by adding the segment connecting
$v^{k,4}$ and $v^{i, j}$ to divide the central
polygon $\Diamond$ into four simplices, where $\{i,
j, k \}= \{1,2,3 \}$. Each $X_{\Xi_k}$ is a
crepant resolution of $X_\Xi(= S_{G})$, and one 
has the refinement relation  of toric varieties
: $
\Xi \prec \Xi_k \prec \Xi^* $ for $k=1, 2, 3$.
The polyhedral decomposition in the
central core  
$\Diamond$ appeared in the refinements 
is indicated by the following relation, 
$$
\Diamond \prec \Diamond_k \prec \Diamond^* \ , \ \
\ \ k=1, 2, 3 ,
$$
whose pictorial realization is shown in
Fig. \ref{fig:hco}. 
The connection between these three smooth $4$-folds corresponding to 
these different
$\Diamond_k$s can be regarded as the ``flop" relation
of 4-folds,  an analogy to the similar procedure
in birational
geometry of 3-folds \cite{Mo}. Each one is a
``small"\footnote{Here the ``smallness" for a
resolution  means one with the exceptional locus
of codimension $\geq 2$.} resolution of the
4-dimensional  isolated singularity with the defining equation  (\ref{4sing}). Hence we have shown
the following result.
\begin{thm} \label{th:flip}
For $G= A_1(4)$, there are crepant resolutions of
$S_G$ obtained by blowing down the divisor $E$ of ${\rm
Hilb}^G(\CZ^4)$ along $(\ref{eq:proj})$ in Theorem
$\ref{th:A1(4)}$. Any two such resolutions differ
by  a ``flop" of $4$-folds. 
\end{thm}

\section{ G-Hilbert Scheme, Crepant Resolution of ${\bf \CZ^4/A_r(4)}$ }
In this section, we give a complete
proof of a general result as in Theorem \ref{th:flip}, but on the group
$A_r(4)$ for all $r$. 

\begin{figure}[ht]
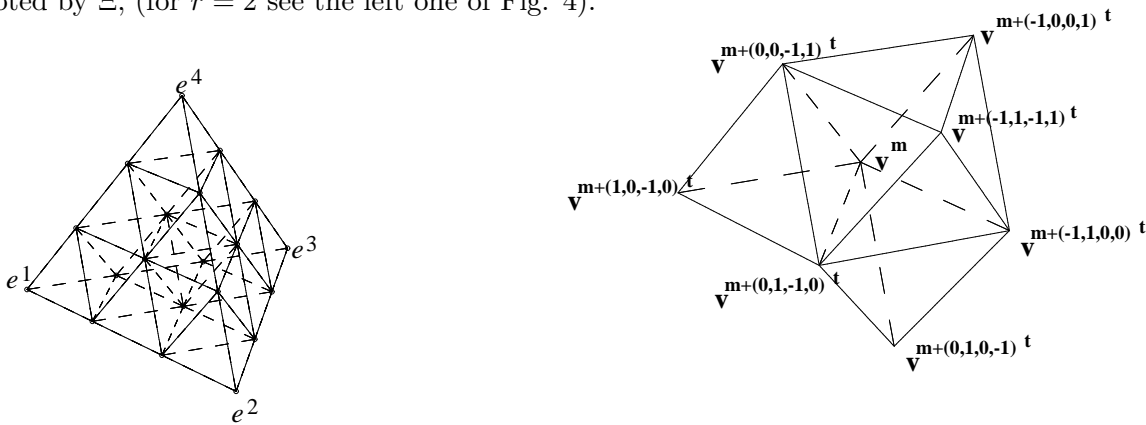

\begin{center}
\mbox{\psfig{figure=tetradec.ps}}
\hskip 2in
\mbox{\psfig{figure=tetraloc.ps}}
\label{fig:tetradl}
\caption{The polytope decomposition $\Xi$ of
$\Delta$ for $r=2$ and local figure of
$\Xi$.} 
\end{center}
\end{figure}

\begin{thm}\label{th:A(4)}
For $G= A_r(4)$, 
the $G$-Hilbert
scheme $\hl^G(\CZ^4)$ is a 
non-singular toric variety  with the
canonical bundle, $
\omega  = {\cal
O}_{\hl^G(\CZ^4)}  ( 
\sum_{k=1}^m E_k ) $ with $ m =
\frac{r(r+1)(r+2)}{6}$,  
where $E_k$s are disjoint smooth
exceptional divisors in
$\hl^G(\CZ^4)$, each of which satisfies the
conditions
(\ref{eq:E}) (\ref{eq:f(-1)}). By blowing down $E_k$ 
to $\PZ^1 \times \PZ^1$ via a
projection (\ref{eq:proj}) for each $k$, it gives rise
to a toric  crepant resolution
$\widehat{S}_G$ of
$S_G$ with  $
\chi ( \widehat{S}_G ) = |A_r(4) |
= (r+1)^3$. 
Furthermore, any two such
$\widehat{S}_G$s differ by a sequence
of flops.
\end{thm}
{\it Proof.}  First we  define the simplicial 
decomposition $\Xi^*$ of (\req(tri)) for
$n=4$, and then we will show  that the
toric variety $X_{\Xi^*}$ is isomorphic to  $\hl^G(\CZ^4)$. We shall denote an element of
$N\cap
\Delta$  by
$$
\vv^{m}(=\vv^{(m_1,..,m_4)^t}) :=
{{m_1e^1+m_2e^2+m_3e^3+m_4e^4}\over {r+1}}  \ ,
\ \ 0 \leq m_i \leq  r+1 \ , \ \sum_{i=1}^4 m_i
= r+1 \ .
$$
For each $\vv^{m} \in N\cap \Delta$, there are
four hyperplanes passing through $\vv^{m}$,  and
parallel to one of the four facets of $\Delta$.
The collection of all such hyperplanes gives rise
to a polytope decomposition of $\Delta$, denoted
by $\Xi$, (for $r=2$ see the left one of Fig. \ref{fig:tetradl}).

Now we examine
the polytope structure of 
$\Xi$. We have $\Xi(0) = N \cap
\Delta$. For each $\vv^{m} \in \Xi(0)$, there
are at most twelve segments in $\Xi(1)$ containing
$\vv^{m}$, and they are given by $\langle \vv^{m},
\vv^{m(i,j)} \rangle$  for 
$i\ne j$,
$1\le i,j\le 4$, 
where $m(i,j) : =  m + e^i-e^j$. 
For a given $\langle \vv^{m}, \vv^{m(i,j)}
\rangle$, the hyperplane passing $\vv^m$ in
$\RZ^4$ with the normal vector $e^i-e^j$ 
separates $\Delta$ into two polytopes $\Delta^\prime$s,  
(one of which could possibly be the empty
set). We are going to discuss those elements in
$\Xi$ containing
$\vv^m$ and lying in a  non-empty polytope of
these two divided ones. For easier description
of our conclusion, also for the  simplicity of 
notions, we shall work on a special model case, say
$i=2, j=3$, and the non-empty polytope
$\Delta^\prime$ consisting of those elements in
$\Delta$ with non-negative inner-product to
$e^2-e^3$, ( no difficulties for a similar
discussion  will arise on other cases except for a
suitable change of indices). The elements in
$\Xi(3)$ contained in $\Delta^\prime$ with
$\vv^m$ as one of its vertices are the following
ones:
\bea(l)
\Delta_u := \langle \vv^m, \vv^{m(2,3)},
\vv^{m(1,3)},\vv^{m(4,3)}
\rangle , \ \ \ \ \ \ \
\Delta_d := \langle \vv^m , 
\vv^{m(2,3)},\vv^{m(2,1)},
\vv^{m(2,4)} \rangle , \\
\Diamond_+ :=
\langle \vv^m , \vv^{m(2,3)},\vv^{m (4,3)},
\vv^{m(2,1)} ,\vv^{m (4,1)} ,
\vv^{m+(-1,1,-1,1)^t}
\rangle ,
\\
\Diamond_- :=
\langle \vv^m , \vv^{m(2, 3)}, \vv^{m(1, 3)},
\vv^{m(2,4)} ,\vv^{m(1,4)} ,
\vv^{m+(1,1,-1,-1)^t}
\rangle \ .
\elea(udD)
Note that $\Diamond_\pm$ are similar by
interchanging $e^3$ and $e^4$, ( for the
configuration of $\Delta_u, \Delta_d, \Diamond_+$,
see the right one of Fig. \ref{fig:tetradl}). Both $\Delta_u, \Delta_d$
are  3-simplices with their vertices forming an integral
basis of $N$, and one facet
of each of these 3-simplices is parallel to that of $\Delta$.
The toric data of $\Delta_u, \Delta_d$ give
rise to the smooth affine open subsets of
$X_{\Xi}$.   The polytope $\Diamond_+$  ($\Diamond_-$) is an octahedron with the
center $c= 
\vv^m+ \frac{e^2+e^4-e^1-e^3}{2(r+1)}$ ($c=\vv^m+
\frac{e^1+e^2-e^3-e^4}{2(r+1)}$ respectively).
We shall mark the octahedron by its center
$c$, and denote it by $\Diamond^c$. The affine
open subset of $X_{\Xi}$ with the
toric data $\Diamond^c$ is smooth except one
isolated singular point $x_{\Diamond^c}$, an 
0-dimensional toric orbit of the affine
toric variety. Hence, one can
conclude that $\Xi(3)$ consists of three type  of elements: $\Delta_u, \Delta_d$ or $
\Diamond^c$. The toric variety
$X_{\Xi}$ is smooth except the finite number 
isolated singularities, $x_{\Diamond^c}$s. The structure of $X_{\Xi}$ near a
singular element  $x_{\Diamond^c}$ can be determined  in the
following manner. For a given $\Diamond^c$, one
can construct a tetrahedron $\Delta^c$ inside
$\Delta$ with the core
$\Diamond^c$ adjacent to four elements 
$\Delta^c_j \ (1\leq j \leq 4)$ in
$\Xi(3)$ of type $\Delta_u$ or $\Delta_d$, 
$$
\Delta^c = \Diamond^c \cup
\bigcup_{j=1}^4 \Delta^c_j
\subseteq \Delta \ , 
$$
such that $\Diamond^c \cap \Delta^c_j \ (1\leq j
\leq 4)$ are four facets of $\Diamond^c $, two
of which intersect
only at one common vertex,  ( 
there could have two possible ways of forming such
$\Delta^c$  with the same core
$\Diamond^c$).  Consider
the rational simplicial decomposition $\Xi^*$ of
$\Delta$, which is a refinement of $\Xi$ 
by adding $c$ as a vertex with the barycentric
simplicial decomposition $\Diamond^c$ for all
$c$. In fact,  the octahedron $\Diamond^c$ is
decomposed into the following eight 4-simplices
of $\Xi^*$:
\bea(ll)
{C_1}^c:=& \langle  c, \ \ 
c+{{e^1+e^2-e^3-e^4}\over 2(r+1)} , \ \
c+{{e^1-e^2+e^3-e^4 }\over2(r+1)}, \ \
c+{{e^1-e^2-e^3+e^4 }\over 2(r+1) }
\rangle , \\ 
{C_2}^c:=& \langle
c+{{e^1+e^2-e^3-e^4}
\over 2(r+1)}, \ \  c ,
\ \ c+{{-e^1+e^2+e^3-e^4}\over 2(r+1)}
, \ \ c+{{-e^1+e^2-e^3+e^4}\over
2(r+1)} \rangle , \\
{C_3}^c:=& \langle
c+{{e^1-e^2+e^3-e^4 }\over 2(r+1)}, \ \
c+{{-e^1+e^2+e^3-e^4}\over 2(r+1)}, \
\  c,
\ \ c+{{-e^1-e^2+e^3+e^4}\over 2(r+1)}
\rangle , \\
{C_4}^c:=& \langle
c+{{e^1-e^2-e^3+e^4}\over 2(r+1)}, \ \
c+{{-e^1+e^2-e^3+e^4}\over 2(r+1)}, \ \
c+{{-e^1-e^2+e^3+e^4}\over 2(r+1)}, \
\ c \rangle ,
\\
{C_1^\prime}^c:=& \langle  c,
\ \ c+{{-e^1-e^2+ e^3 +e^4}\over
2(r+1)},
\
\ c+{{-e^1+ e^2 -e^3+ e^4}\over
2(r+1)},
\ \ c+{{-e^1+ e^2+ e^3- e^4}\over
2(r+1)}
\rangle , \\
{C_2^\prime}^c:=& \langle
c+{{-e^1-e^2+e^3+e^4}\over 2(r+1)}, \ \
c ,
\
\ c+{{e^1-e^2-e^3+e^4}\over 2(r+1)} , \
\ c+{{e^1-e^2+e^3-e^4}\over 2(r+1)}
\rangle , \\
{C_3^\prime}^c:=& \langle
c+{{-e^1+e^2-e^3+e^4 }\over 2(r+1)}, \
\ c+{{e^1-e^2-e^3+e^4 }\over 2(r+1)}, \
\  c, \ \ 
c+{{e^1+e^2-e^3-e^4 }\over 2(r+1)}
\rangle , \\ 
{C_4^\prime}^c:=& \langle
c+{{-e^1+e^2+e^3-e^4 }\over 2(r+1)}, \
\ c+{{e^1-e^2+e^3-e^4}\over 2(r+1)}, \
\ c+{{e^1+e^2-e^3-e^4}\over 2(r+1)}, c
\rangle . 
\elea(Cc)
All vertices appeared in the
above simplices are
elements in
$N \cap \Delta$ except $c$ , while  $2c \in N$.
(see Fig. 5)

\begin{figure}[ht]
\begin{center}
\mbox{\psfig{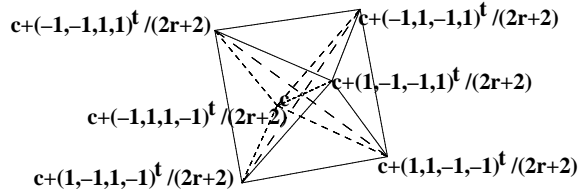}}
\label{fig:hilbst}
\caption{Local figure of the decomposition of the 
octahedron in the right one of Fig. \ref{fig:tetradl} by  adding $c$.} 
\end{center}
\end{figure}
One can determine the singularity structure
of the variety $X_{\Xi}$ near
$x_{\Diamond^c}$ by examining the toric orbits
associated to $\Delta^c$.
The toric data in $\RZ^4$ for the lattice
$N$ and the cone  generated by
$\Delta^c$ are isomorphic to the toric data of
the lattice for the group
$A_1(4)$ with the first quadrant cone in Lemma 
\ref{lem:SingXi}. Hence as toric varieties,
the structure of $X_{\Xi}$ near the singularity
$x_{\Diamond^c}$ is the same as that for $A_1(4)$.
One can apply the result of Theorem
\ref{th:A(4)} to describe the local
structure of
$X_{\Xi^*}$ over the singular point
$x_{\Diamond^c}$ of
$X_{\Xi}$. Hence one concludes that  $X_{\Xi^*}$
is a smooth toric variety with the canonical
bundle,  $
\omega_{X_{\Xi^*}} = {\cal O}_{X_{\Xi^*}} ( 
\sum_{\Diamond^c \in \Xi(4)} E_c )$,  
where $E_c$ is the toric divisor associated to
the vertex $c$ in
$X_{\Xi^*}$, and it satisfies the properties
$(\ref{eq:E}) (\ref{eq:f(-1)})$. By
(\ref{tori}) and the structure of
$E_c$, one obtains the
desired crepant resolutions $\widehat{S}_{A_r(4)}$
by  blowing-down each $E_c$ to 
$\PZ^1 \times \PZ^1$ as in Theorem \ref{th:flip}. and different crepant resolutions are connected by flop relation. 
It remains to show 
$X_{\Xi^*}
\simeq \hl^G(\CZ^4)$, and the total number of $\Diamond^c$s is equal to
$\frac{r(r+1)(r+2)}{6}$. As in the proof of
Theorem \ref{th:A1(4)}, we first construct a
regular morphism $\lambda$ from $X_{\Xi^*}$ to
$\hl^G(\CZ^4)$ by examining 
$I(y)$ for $y \in X_{\Xi^*}$ in terms of toric
coordinates. For $\RT \in \Xi^*(3)$,
we denote $x_\RT : = {\rm orb}(\RT) \in
X_{\Xi^*}$. For the simplicity  of notions, we
again work on some special 3-simplices as 
the model cases, whose  argument can equally be
applied to all  elements in
$\Xi^*(3)$. We consider the 3-simplices of
$X_{\Xi^*}$ contained in the first three
polytopes in (\req(udD)), and they are: 
$\Delta_u, \Delta_d$ of (\req(udD)) and the
eight simplices of (\req(Cc)) with $c= 
\vv^m+ \frac{e^2+e^4-e^1-e^3}{2(r+1)}$. The affine
toric coordinates for $X_{\Xi^*}$ are determined
by the integral basis of $M$ in the 
simplicial cone dual to the one in $N$ generated by the
corresponding 3-simplex.  By computation, the
affine coordinate systems corresponding to these
3-simplices are as follows:
\begin{eqnarray*}
\Delta_u : & ( V_1^{(m_1)} ,  
V_2^{(m_2)} ,  V_3^{(m_3-1)} , 
V_4^{(m_4)} ) , & \ \  V_i^{(l)}:= \ {
{{Z_i}^{r+1-l}} \over { ( Z_1..\breve{Z}_i..Z_4
)^l} } , \\
\Delta_d : & 
( V_1^{\prime (m_1)} ,  
V_2^{\prime (m_2+1)} ,  V_3^{\prime (m_3)} , 
V_4^{\prime (m_4)} ) , & \ \ 
V_i^{\prime (l)}: = \ { {( Z_1..\breve{Z}_i..Z_4
)^l}   \over {  {{Z_i}^{r+1-l}} }} \ , \\
{C_i}^c : & (
U_{i,1}^{(c)}, \ U_{i,2}^{(c)}, \ U_{i,3}^{(c)}, \
U_{i,4}^{(c)} ) , & U_{i,i}^{(c)}:= 
{ {( Z_jZ_jZ_k)^{(r+1)c_i+\frac{1}{2}}}   \over { 
{{Z_i}^{(r+1)(1-c_i)-\frac{1}{2}}} }} , \ 
U_{i,j}^{(c)}:= 
{ {( Z_iZ_j)^{(r+1)(1-c_i-c_j)}}   \over { 
{{Z_kZ_s}^{(r+1)(c_i+c_j)}} }}, \\
C_i^{\prime c}:& (
U_{1,i}^{\prime (c)}, \ U_{2,i}^{\prime (c)}, \
U_{3,i}^{\prime (c)}, \ U_{4,i}^{\prime (c)} ) , &
U_{i,i}^{\prime (c)}:=
{{{Z_i}^{(r+1)(1-c_i)+{1\over 2}}}
\over{(Z_jZ_kZ_s)^{(r+1)c_i-\frac{1}{2}}}}, \ 
U_{i, j}^{\prime (c)}:={{
(Z_kZ_s)^{(r+1)(1-c_k-c_s)}}\over
{(Z_iZ_j)^{(r+1)(c_k+c_s)}}}. 
\end{eqnarray*} 
Here the indices $i,j, k, s$ indicate the
four 3 by permuting $1,2,3,4$ ,
and we shall adopt this convention for the rest of
this proof if no confusion will arise. Define
the following eigen-polynomials of
$G$ for
$\beta \in \CZ$ and integers $l$ with $ 0 \leq l \leq (r+1)$,
$$
F_i^{(l)}(\beta)= 
{Z_i}^l-\beta(Z_jZ_kZ_s)^{(r+1)-l},  
G_{i,j}^{(l)}(\beta)= 
{(Z_iZ_j)}^l-\beta (Z_kZ_s)^{(r+1)-l}, 
H_i^{(l)}(\beta)= 
{(Z_jZ_jZ_s)}^l-\beta{Z_i}^{(r+1)-l}.
$$ 
Let $y$ be an element of $X_{\Xi^*}$. For $y$
near $x_{\Delta_u}$ with coordinates $
(V_1^{(m_1)}, V_2^{(m_2)}, V_3^{(m_3-1)}, V_4^{(m_4)}) =(v_1, v_2, v_3, v_4)$, 
the ideal $I(y)$ has the generators,
\bea(l)
F_1^{(r+1-m_1)}(v_1), \ 
F_2^{(r+1-m_2)}(v_2), \
F_3^{(r+2-m_3)}(v_3), \
F_4^{(r+1-m_4)}(v_4), \
G_{1,2}^{(m_3+m_4)}(v_1v_2), \\
G_{1,3}^{(m_2+m_4+1)}(v_1v_3), \
G_{1,4}^{(m_2+m_3)}(v_1v_4),
\ 
G_{2,3}^{(m_1+m_4+1)}(v_2v_3),
\
G_{2,4}^{(m_1+m_3)}(v_2v_4),
\ G_{3,4}^{(m_1+m_2+1)}(v_3v_4),
\\
H_1^{(m_1+1)}(v_2v_3v_4),
\ H_2^{(m_2+1)}(v_1v_3v_4),
\ H_3^{(m_3)}(v_1v_2v_4),
\ H_4^{(m_4+1)}(v_1v_2v_3), \
Z_1Z_2Z_3Z_4-v_1v_2v_3v_4.
\elea(Delu)
For $y$ near $x_{\Delta_d}$ with coordinates 
$
(V_1^{\prime (m_1)},   V_2^{\prime (m_2+1)},  
V_3^{\prime (m_3)} ,  V_4^{\prime
(m_4)})=(v_1^\prime , v_2^\prime , v_3^\prime, v_4^\prime)$, 
$I(y)$ has the generators:
\bea(l)
F_1^{(r+2-m_1)}(v^\prime_2v^\prime_3v^\prime_4),
\
F_2^{(r+1-m_2)}(v^\prime_1v^\prime_3v^\prime_4)
, \
F_3^{(r+2-m_3)}(v^\prime_1v^\prime_2v^\prime_4)
,
\ 
F_4^{(r+2-m_4)}(v^\prime_1v^\prime_2v^\prime_3),\\
G_{1,2}^{(m_3+m_4)}(v^\prime_3v^\prime_4),
\
G_{1,3}^{(m_2+m_4+1)}(v^\prime_2v^\prime_4),
\ G_{1,4}^{(m_2+m_3+1)}(v^\prime_2v^\prime_3),
\
G_{2,3}^{(m_1+m_4)}(v^\prime_1v^\prime_4),
\
G_{2,4}^{(m_1+m_3)}(v^\prime_1v^\prime_3),
\\ G_{3,4}^{(m_1+m_2+1)}(v^\prime_1v^\prime_2), \
H_1^{(m_1)}(v^\prime_1),
\
H_2^{(m_2+1)}(v^\prime_2),
\ H_3^{(m_3)}(v^\prime_3), \
H_4^{(m_4)}(v^\prime_4), \
Z_1Z_2Z_3Z_4-v^\prime_1v^\prime_2v^\prime_3v^\prime_4.
\elea(Deld)
For $y$ near $x_{C^c_i}$ with coordinates 
$(U_{il}^{(c)}=u_l)_{1 \leq l \leq 4} $,
$I(y)$ has the generators:
\bea(l)
F_i^{((r+1)(1-c_i)+{1\over
2})}(u_1u_2u_3u_4),
\ F_j^{((r+1)(1-c_j)+{1\over
2})}(u_iu_j) ,
\ F_k^{((r+1)(1-c_k)+{1\over
2})}(u_iu_k), \\ 
F_s^{((r+1)(1-c_s)+{1\over 2})}(u_iu_s),
\
G_{i,j}^{(r+1)(c_k+c_s)}(u_j), \ 
G_{i,k}^{(r+1)(c_j+c_s
)}(u_k), \ G_{i,s}^{(r+1)(c_j+c_k)}(u_s),
\\
G_{j,k}^{(r+1)(c_i+c_s)+1}({u_i}^2u_ju_k),
\
G_{j,s}^{(r+1)(c_i+c_k)+1}({u_i}^2u_ju_s),
\ G_{k,s}^{(r+1)(c_i+c_j)+1}({u_i}^2
u_ku_s), \\
H_i^{((r+1)c_i+{1\over
2})}(u_i), \ H_j^{((r+1)c_j+{1\over
2})}(u_iu_ku_s), \
H_k^{((r+1)c_k+{1\over
2})}(u_iu_ju_s), \
H_s^{((r+1)c_s+{1\over 2})}(u_iu_ju_k),
\\
Z_1Z_2Z_3Z_4-{u_i}^2u_ju_ku_s.
\elea(Cci)
For $y$ near $x_{C^{\prime c}_i}$ with the
coordinates 
$(U_{il}^{\prime (c)}=u^\prime_l)_{1 \leq l \leq
4}
$, $I(y)$ has the generators:
\bea(l)
F_i^{((r+1)(1-c_i)+{1\over
2})}(u^\prime_i)\ , F_j^{((r+1)(1-c_)+{1\over
2})}(u^\prime_iu^\prime_ku^\prime_s),
\ F_k^{((r+1)(1-c_k)+{1\over
2})}(u^\prime_iu^\prime_ju^\prime_s), \ 
F_s^{((r+1)(1-c_s)+{1\over
2})}(u^\prime_iu^\prime_ju^\prime_k),\\
G_{i,j}^{(r+1)(c_k+c_s)+1}({u^\prime_i}^2
u^\prime_ku^\prime_s),
\
G_{i,k}^{(r+1)(c_j+c_s)+1}({u^\prime_i}^2
u^\prime_ju^\prime_s),
\
G_{i,s}^{(r+1)(c_j+c_k)+1}({u^\prime_i}^2
u^\prime_j
u^\prime_k),\\ 
G_{j,k}^{(r+1)(c_i+c_s)}(u^\prime_s),
\ G_{j,s}^{(r+1)(c_i+c_k )}(u^\prime_k), \
G_{k,s}^{(r+1)(c_i+c_j)}(u^\prime_j), \\
H_i^{((r+1)c_i+{1\over
2})}(u^\prime_1u^\prime_2u^\prime_3u^\prime_4),
\
H_j^{((r+1)c_j+{1\over
2})}(u^\prime_iu^\prime_j), \
H_k^{((r+1)c_k+{1\over
2})}(u^\prime_iu^\prime_k)\mbox{,
}H_s^{(c_s+{1\over 2})}(u^\prime_iu^\prime_s),
\\
Z_1Z_2Z_3Z_4-{u^\prime_i}^2u^\prime_j
u^\prime_ku^\prime_s.
\elea(Ccip)
The centers of the above affine charts have the monomial ideals , say the one
near $x_{\Delta_u}$, $I (x_{\Delta_u})$ is obtained by setting $v_l=0$ in (\req(Delu)), hence an monomial
ideal. There are exactly
$(r+1)^3$ monomials not in $I (x_{\Delta_u})$,
i.e., 
$ | I (x_{\Delta_u})^\dagger | = (r+1)^3 $. For $y$
near $x_{\Delta_u}$, by using (\req(Delu)) and employing
the Gr\"{o}bner basis techniques and the toric
data, one obtains the colength of
$I(y)$ in $\CZ[Z]$ satisfying the relation,
${\rm colength}(I(y))
\le{\rm colength}(I(x_{\Delta_u}))=(r+1)^3$; 
this implies ${\rm colength}(I(y))=(r+1)^3$. By which
it determines an element $\lambda(y) \in {\rm
Hilb}^G(\CZ^4)$. One can
also show the colength of $I(y)$ equal to $(r+1)^3$
for $y$ in other affine charts using
(\req(Deld)) (\req(Cci)) (\req(Ccip)). The
same conclusion holds for $y$ in any affine
coordinate neighborhood centered at $x_\RT$ for
$\RT \in
\Xi^*(3)$, and one obtains an element $\lambda(y)$
in 
$\hl^G(\CZ^4)$, by which the
morphism 
$\lambda: X_{\Xi^*} \longrightarrow {\rm
Hilb}^G(\CZ^4)$ is defined. 

Now we are going to show that
$\lambda$ is an isomorphism. For  $n \in \ZZ$,
we denote $\underline{n}$ the
unique integer satisfying the relation, 
$$
n \equiv \underline{n} \pmod{r+1} \ , \ \ 
0 \leq \underline{n} \leq r  \ .
$$
We first determine the
$G$-invariant monomial ideals $J_0$ in 
$\hl^{A_r(4)}(\CZ^4)$. For a such
$J_0$,  the set $\eM := W_0 \setminus ( W_0 \cap
J_0 ) $ forms a basis of a $G$-regular
representation space.
Denote
$l_i$ the smallest integer with 
$Z_i^{l_i} \in J_0$; $l_{ij}$ the smallest
one with $(Z_i Z_j)^{l_{ij}} \in J_0$ for $i
\ne j$,  and so on. 
By $1
\not\in J_0$, and
$1\sim {Z_i}^{r+1}\sim Z_1Z_2Z_3Z_4$, we have 
${Z_i}^{r+1},Z_1Z_2Z_3Z_4 \in J_0$, i.e. $I (o)
\subset J_0$, and the following relations hold,
$$
1 \leq l_{i j k} \leq l_{i j} \leq l_i \leq
r+1 \ .
$$ 
By $J_0^\bot
\subset I(o)^\bot$, and 
(\req(Ioorg)) for the description of the
$G$-eigenspace of $I(o)^\bot$, 
$(Z_jZ_kZ_s)^{r+1-l_i}$ is the only monomial $u
\in I(o)^\dagger$ with $u \sim Z_i^{l_i}$, which
implies $
(Z_jZ_kZ_s)^{r+1-l_i}\in \eM $ and $ 
(Z_jZ_kZ_s)^{r+2-l_i}\in J_0 $, hence $l_{jks}= r+2-l_i$. By a similar
argument, one has $ l_{ks} = r+2-l_{ij}$. Hence we have
\be
l_i + l_{jks}=  l_{ij} + l_{k s} =
r+2 \ .
\ele(ijks)
We claim that $J_0$ is the ideal with 
generators given by 
\be
J_0= \langle {Z_i}^{l_i}, \ (Z_iZ_j)^{l_{ij}},
\ (Z_iZ_jZ_k)^{l_{i j k}}, \ Z_1Z_2Z_3Z_4 \
\ | \ i, j , k \rangle \ .
\ele(J0)
(Note that $i, j, k$ are distinct numbers among
$1,2,3,4$ as before).  Let $J_0'$ the ideal in the
right hand side of (\req(J0)
). Then 
$I(o)
\subset J_0' 
\subset J_0$. Suppose 
$J_0'  \ne J_0$, equivalently $J_0 \cap J_0'^\dagger
\ne \emptyset$. For the convenience of
notations but without loss of generality, we may
assume
$Z_2^{i_2}Z_3^{i_3}Z_4^{i_4} \in J_0 \cap
J_0'^\dagger$ for $i_2 \leq i_3 \leq i_4$. Hence 
$i_2 < l_{234} ,   i_3 < l_{34} ,  i_4 <
l_4$,
which implies
$p_1 \  (: =
Z_2^{l_{234}-1}Z_3^{l_{34}-1}Z_4^{l_4-1}) 
\in J_0 \cap I(0)^\dagger $.   
By (\req(Ioorg)), the rest of monomials $p$ in 
$I(o)^\dagger$ with $p \sim p_1$ are
given by
$$
\begin{array}{c}
p_2:= Z_1^{\underline{r+2-l_{234}}}
Z_3^{l_{34}-l_{234}}Z_4^{l_4-l_{234}}
\ , \ \ 
p_3:= Z_1^{\underline{r+2-l_{34}}}
Z_2^{\underline{r+1+l_{234}-l_{34}}}Z_4^{l_4-l_{34}}
,
\\
p_3:= Z_1^{\underline{r+2-l_4}}Z_2^{
\underline{r+1+l_{234}-l_4}}Z_3^{\underline{
r+1+l_{34}-l_4}},
\end{array}
$$
among which exactly only one belongs to
$J_0^\dagger$. We have $p_1=p_2$ when $l_{234}=1$.
If $l_{234} > 1$, by (\req(ijks)) we have 
$\underline{r+2-l_{234}}= l_1$.
Therefore $p_2 \in J_0$. When
$l_{234}=l_{34}$, we have $p_2= p_3$.
When $l_{234}<l_{34}$, $p_3 = 
(Z_1Z_2)^{l_{1,2}}
Z_2^{l_{123}}Z_4^{l_4-l_{34}}$ by (\req(ijks)),
hence $p_3 \in J_0$. Similarly, $p_3= p_4$ when
$l_{34}=l_4$. If $l_{34}<l_4$,
$u_4=(Z_1Z_2Z_3)^{l_{123}}Z_2^{l_{234}}
Z_3^{l_{34}}$, hence $ p_4 \in J_0$. Therefore
$p_i \in J_0$ for $1 \leq i \leq 4$, a   
contradiction to their relations with $J_0^\dagger$. We
are going to show the following relations hold for
$i \ne j$,
\be
r+1 \leq l_i + l_j - l_{ij} \leq r+2 \ .
\ele(ij)
Consider the element $w \ (:=
Z_i^{l_i}Z_j^{l_{ij}-1}Z_k^{l_{ijk}-1})$ in
$J_0$. Among the following monomials
$G$-equivalent to $w$,
$$
\begin{array}{c}
w_1 =
Z_i^{\underline{l_i-l_{ijk}+1}}
Z_j^{l_{ij}-l_{ijk}}
Z_s^{\underline{r+2-l_{ijk}}} \ , \ \
w_2 = Z_i^{\underline{l_i-l_{ij}+1}}
Z_k^{\underline{r+1-l_{ij}+l_{ijk}}}
Z_s^{\underline{r+2-l_{ij}}} \ , \\
w_3 = Z_j^{r-l_i+l_{ij}} Z_k^{r-l_i+l_{ijk}}
Z_s^{r+1-l_i} ,
\end{array}
$$ 
there exists exactly one in $J_0^\dagger$. It is
easy to see that 
$w_1 = Z_i^{\underline{l_i-l_{ijk}+1}}
Z_j^{l_{ij}-l_{ijk}}
Z_s^{l_s} \in J_0$ unless $l_{ijk}=1$,
in which case
$w_1=w \in J_0$ if $l_i <r+1$, and $w_1=w_3$ if
$l_i=r+1$. 
We have
$w_1=w_2$ if $l_{ij} = l_{ijk}$. When $l_{ij} >
l_{ijk}$,  
$w_2 =  Z_i^{\underline{l_i-l_{ij}+1}}
Z_k^{l_{ks}+l_{ijk}-1}
Z_s^{l_{ks}} \in J_0$. Therefore $w_3$ is the
element of
$J_0^\dagger$  $G$-equivalent to $w$, which by
the expression of the power of $Z_j$, implies  
$$
r+1 \leq l_i+l_j -l_{ij} \ .
$$
As a consequence of the above inequality,  we
have $l_j= r+1$ and
$l_i+l_j-l_{ij}= r+1$ when 
$l_{ij}=l_i$, in particular 
(\req(ij)) holds. Hence we may assume $l_{ij}<
l_i$, in which case $h :=
Z_i^{l_i-1}Z_j^{l_{ij}}Z_k^{l_{ijk}-1} \in J_0$.
Among the following monomials
$G$-equivalent to $h$,
$$
\begin{array}{c}
h_1 =
Z_i^{l_i-l_{ijk}}
Z_j^{l_{ij}-l_{ijk}+1}
Z_s^{\underline{r+2-l_{ijk}}} \ , \ \
h_2 = Z_i^{l_i-l_{ij}-1} Z_k^{r-l_i+l_{ijk}}
Z_s^{r+1-l_{ij}} , \\
h_3 = Z_j^{\underline{r+2-l_i+l_{ij}}}
Z_k^{r+1-l_{j}+l_{ijk}}
Z_s^{r+2-l_{i}} \ , 
\end{array}
$$ 
there exists exactly one  in $J_0^\dagger$.
We have
$h_1=h \in J_0$ if $l_{ijk}=1$. When $l_{ijk}>1$,
$ h_1 =
Z_i^{l_i-l_{ijk}}
Z_j^{l_{ij}-l_{ijk}+1}
Z_s^{l_s}$, and $h_1 \in J_0$. One has
$h_3 = Z_j^{l_{ij}} Z_k^{l_{ijk}-1}
(Z_jZ_kZ_s)^{l_{jks}} \in J_0 $ unless
$l_i=l_{ij}+1$, in which case, $h_3=h_2$.
Therefore we have $h_2 \in J_0^\dagger$, which
implies $l_i -l_{ij} -1 \leq l_{iks}-1$, hence 
$l_i + l_j -l_{ij} \leq r+2 $ by (\req(ijks)
).
Therefore we obtain the relation (\req(ij)). 
With $(i, j)=(1,2), (3, 4)$ in (\req(ij)) 
(\req(ijks)), we have $
3r+4 \leq \sum_{j=1}^4 l_j \leq 3r+6 $.  
Using (\req(ijks)), one obtains the all possible
cases of $l_i+l_j-l_{ij}$ for a given value of
$\sum_{j=1}^4l_j$; consequently, 
all the values of  $l_I$s are determined by $l_i$s. 
By comparing the polynomials at the
origin in (\req(Delu)) (\req(Deld)) (\req(Cci))
(\req(Ccip)), $J_0= I(x_\RT)$ for $\RT \in
\Xi^*(3)$ by the following relations: 
\bea(lll)
\Delta_u : & \sum_{j=1}^4 l_j = 3r + 4 , & l_{i
j} = l_i + l_j -r-1 ; \\
\Delta_d : & \sum_{j=1}^4 l_j = 3r + 6 , & l_{i
j} = l_i + l_j -r-2 ; \\
C_i^c : & \sum_{j=1}^4 l_j = 3r + 5 , & l_{ij}=
l_i+l_j-r-2 \ , \ l_{ks} = l_k+l_s-r-1 ; \\
C_i^{\prime c} : & \sum_{j=1}^4 l_j = 3r + 5 , &
l_{ij}= l_i+l_j-r-1 \ , \ l_{ks} = l_k+l_s-r-2 \ ,
\elea(lsrel)
where the indices in toric data are connected
to the
$l_i$s by the following relations,
$$
\begin{array}{ll}
\Delta_u :  & 
l_3 = r+2-m_3, \ l_j =r+1-m_j, \ ( j \ne
3 ) \ , \\ 
\Delta_d :  & 
l_2 = r+1-m_2, \ l_j =r+2-m_j, \ ( j \ne 2 ) \ , 
\\
C_i^c , C_i^{\prime c} :  & l_j = (r+1)(1-c_j)+
\frac{1}{2} \ , \ \ \  c = \frac{1 }{2r+2}
\sum_{j=1}^4 (2 r+3 -2l_j ) e^j \ .
\end{array}
$$
With $l_i':= r+1-l_i$,
$l_i'$s are 4 positive integers satisfying the
equation $\sum_{i=1}^4 l_i' = L'$ with $L'=r,
r-1, r-2$. The number of solutions of $l_i'$s is equal to   
$\left(\matrix{{L'+3}\cr 3}\right)$. Hence one
obtains the following  numbers of
$\RT \in \Xi^*(3)$ for the toric data in
(\req(udD)) (\req(Cc)) using the relation with  
$l_i$s:
\be
\# \{ \Delta_u \} = \frac{(r+1)(r+2)(r+3)}{6}, \
\ 
\#\{\Delta_d \} = \frac{(r-1)r(r+1)}{6}, \ \ 
\#\{ c \} = 
\frac{r(r+1)(r+2)}{6} . 
\ele(numbty)

Let $J$ be a $G$-invariant ideal representing an
element in $\hl^G(\CZ^4)$. With the 
Gr\"{o}bner basis argument as in Theorem
\ref{th:A1(4)}, there is a monomial ideal $J_0$ in
$\hl^G(\CZ^4)$ such that $J_0^\dagger$
gives rise to a basis of $\CZ[Z]/J$ with the
relation (\req(Jgen)). As $J_0 = I (x_\RT)$ for some $\RT
\in
\Xi^*(3)$, which is determined by the integers
$l_i, l_{ij}, l_{ijk}$ with the relations
in (\req(ijks)) (\req(lsrel)), 
this implies that for some 
$\gamma_i, \gamma_{ij}, \gamma_{jks},
\gamma_{1234} \in \CZ$, the polynomials 
$F_i^{(l_i)}(\gamma_i),
G_{ij}^{(l_{ij})}(\gamma_{ij}),
H_i^{(l_{jks})}(\gamma_{jks})$ and
$Z_1Z_2Z_3Z_4-\gamma_{1234}$ are elements of
$J$. From the expressions of  $F_i^{(l)}(\beta), 
G_{i,j}^{(l)}(\beta),
H_i^{(l)}(\beta)$, and using ${\rm
dim}(\CZ[Z]/J)= (r+1)^3$, one can conclude
$$
J = \langle F_i^{(l_i)}(\gamma_i),
G_{ij}^{(l_{ij})}(\gamma_{ij}),
H_i^{(l_{jks})}(\gamma_{jks}),
Z_1Z_2Z_3Z_4- \gamma_{1234}  \rangle_{i, j ,
k , s}
$$
We are going to
determine the relations among the $\gamma_I$s using the relations
(\req(ijks))(\req(lsrel)) and 
according to the type of
$l_i$s.  By 
$$
(\gamma_{1234}-\gamma_{123}\gamma_4) 
{Z_4}^{l_4-1}
= Z_1Z_2Z_3 F_4^{(l_4)}(\gamma_4)-
\gamma_4 H_4^{(l_{123})}(\gamma_{123}) 
-{Z_4}^{l_4-1}(Z_1Z_2Z_3Z_4-\gamma_{1234}) \in
J ,
$$
and ${Z_4}^{l_4-1} \not\in J$, we have
$\gamma_{1234}= \gamma_{123}\gamma_4 $. 

For $J$ with $J_0$ of type $\Delta_u$, by
(\req(lsrel)) we have 
$$
\begin{array}{ll}
(\gamma_{123}-\gamma_{12} \gamma_3) 
{Z_3}^{l_{34}-1}{Z_4}^{l_4-1}=
{(Z_1Z_2)^{l_{123}}
F_3^{(l_3)}(\gamma_3)
+\gamma_3{Z_4}^{l_{124}-1}
G_{12}^{(l_{12})}(\gamma_{12})
-{Z_3}^{l_{34}-1}H_4^{(l_{123})}(\gamma_{123}) ,}
\\
(\gamma_{13} -\gamma_1
\gamma_3){Z_3}^{l_{234}-1}(Z_2Z_4)^{l_{24}-1}
= \gamma_3 (Z_2Z_4)^{l_{124}-1}
F_1^{(l_1)}(\gamma_1) +{Z_1}^{l_{13}}
F_3^{(l_3)}(\gamma_3)
-{Z_3}^{l_3-l_{13}}G_{13}^{(l_{13})}(
\gamma_{13}) , 
\end{array}
$$
which are elements in $J$. By
${Z_3}^{l_{34}-1}{Z_4}^{l_4-1}, 
{Z_3}^{l_{234}-1}(Z_2Z_4)^{l_{24}-1} \in
J_0^\dagger$, we have $\gamma_{123}=
\gamma_1 \gamma_{23}, \gamma_{2 3} =
\gamma_2\gamma_3$. By permuting the indices, one
obtains $\gamma_I= \prod_{i
\in I} \gamma_i$ for a subset $I$ of
$\{1,2,3,4\}$. By (\req(Delu))(\req(numbty)), we
have $J= I(y)$ for $y$ near $x_{\Delta_u}$ 
with the coordinates $v_i =\gamma_i$.

When $J_0$ is of type $\Delta_u$, by
(\req(lsrel)), the
following elements are in $J$,
$$
\begin{array}{l}
(\gamma_{12}\gamma_{134}-\gamma_{1})
{Z_2}^{l_{234}-1}(Z_3Z_4)^{l_{34}-1}
=(Z_3Z_4)^{l_{134}}
F_1^{(l_1)}(\gamma_1)
-\gamma_{134}{Z_2}^{l_{234}-1}
G_{12}^{(l_{12})}(\gamma_{12}) 
-{Z_1}^{l_{12}}
H_2^{(l_{134})}(\gamma_{134}),  \\
(\gamma_{12}-\gamma_{123}\gamma_{124})
{Z_3}^{l_3-1}{Z_4}^{l_{34}-1}= 
-{Z_3}^{l_{123}}
G_{12}^{(l_{12})}(\gamma_{12})
+\gamma_{123}{Z_4}^{l_{34}-1}
H_3^{(l_{124})}(\gamma_{124})
+(Z_1Z_2)^{l_{124}}
H_4^{(l_{123})}(\gamma_{123}) \ .
\end{array}
$$
Therefore $\gamma_1 =
\gamma_{12}\gamma_{134}$ and
$\gamma_{12}=\gamma_{123}\gamma_{124}$. Set
$v_i^\prime= \gamma_{1..\breve{i}..4}$. With 
the same argument, one
obtains $\gamma_I = \prod_{j \in I'} v_j^\prime$
for $I \ne 1234$, where  $I'$ is the complement
set of $I$ in $\{1,2,3,4\}$. Therefore by
(\req(Deld)) (\req(numbty)), 
$J= I(y)$
for $y$ near $x_{\Delta_d}$  
having $v_i^\prime$'s as coordinates.

When $J_0$ is of type $C^c_i$ or $C^{\prime
c}_i$, without loss of generality, we may assume 
$i=1$. In
the case $C^c_1$, the following elements are in
$J$ by (\req(lsrel)),
$$
\begin{array}{l}
(\gamma_{123}-\gamma_{13}\gamma_2)
(Z_1Z_3)^{l_{134}-1}{Z_4}^{l_4-1}=
\gamma_{13}{Z_4}^{l_{24}-1}F_2^{(l_2)}(\gamma_2)
+{Z_2}^{l_{13}-l_{134}+1}G_{13}^{(l_{13})}(
\gamma_{13})
-(Z_1Z_3)^{l_{134}-1}H_4^{(l_{123})}(
\gamma_{123}) , \\
(\gamma_2-\gamma_{12}\gamma_{234})Z_1^{l_{134}-1}
(Z_3Z_4)^{l_{34}-1}=
-(Z_3Z_4)^{l_{234}}
F_2^{(l_2)}(\gamma_2) 
+\gamma_{234}{Z_1}^{l_{134}-1}
G_{12}^{(l_{12})}(\gamma_{12})
+Z_2^{l_{12}}H_1^{(l_{234})}(\gamma_{234}) , \\
(\gamma_1-\gamma_{12}\gamma_{134})Z_2^{l_{234}-1}
(Z_3Z_4)^{l_{34}-1}
=
-(Z_3Z_4)^{l_1-l_{1,2}}
F_1^{(l_1)}(\gamma_1)
+\gamma_{134}{Z_2}^{l_{234}-1}
G_{12}^{(l_{12})}(\gamma_{1,2})
+{Z_1}^{l_{12}}H_2^{(l_{134})}(\gamma_{134}) , \\
(\gamma_{23}-\gamma_2\gamma_3)(Z_2)^{l_{124}-1}
(Z_1Z_4)^{l_{14}-1} =
\gamma_2
(Z_1Z_4)^{l_3-l_{23}} F_3^{(l_3)}(\gamma_3) 
+{Z_3}^{l_{23}}F_2^{(l_2)}(\gamma_2)
-{Z_2}^{l_{124}-1}
G_{23}^{(l_{23})}(\gamma_{23}) \ .
\end{array}
$$
Hence 
$$
\gamma_{123}= \gamma_2 \gamma_{13} \ ,
\gamma_2= \gamma_{234} \gamma_{12} \ , 
\gamma_1=\gamma_{12}\gamma_{134} \ , 
\gamma_{23}=\gamma_2\gamma_3 ,
$$
which are the same relations as $u_I$s in
(\req(Cci)) for $i=1$ under the identification: 
$u_1=\gamma_{234}$, and $u_j=\gamma_{1j}$
for $j\ne 1$.
By permuting the indices, one can show that
all the rest relations in (\req(Cci)) are
satisfied in terms of the $\gamma_I$s. Hence by
(\req(numbty)), $J=I(y)$ for 
$y$ near $x_{C_1}$ with $u_i$s as the
coordinates of $y$.

For $J_0$ is of type $C^{\prime
c}_1$, the following elements are in
$J$ by (\req(lsrel)), 
$$
\begin{array}{l}
(\gamma_{234}-\gamma_{34}\gamma_2)
{Z_1}^{l_1-1}{Z_2}^{l_{12}-1} =
(Z_3Z_4)^{l_{234}}F_2^{(l_2)}(
\gamma_2)+ \gamma_2{Z_1}^{l_1-l_{12}}
G_{34}^{(l_{34})}(\gamma_{34})
-Z_2^{l_{12}-1}
H_1^{(l_{234})}(\gamma_{234}),

\\
(\gamma_2-\gamma_{23}\gamma_{124})
(Z_1Z_4)^{l_{134}-1}{Z_3}^{l_3-1}
=
-{Z_3}^{l_{23}}F_2^{(l_2)}(\gamma_2) 
+{Z_2}^{l_{124}}
G_{23}^{l_{2,3}}(\gamma_{23})
+\gamma_{23}(Z_1Z_4)^{l_{134}-1}
H_3^{(l_{124})}(\gamma_{124}) ,

\\
(\gamma_{124}-\gamma_1\gamma_{24})
{Z_1}^{l_{13}-1}{Z_3}^{l_3-1}
= (Z_2Z_4)^{l_{124}}
F_1^{(l_1)}(\gamma_1)
+\gamma_1{Z_3}^{l_3-l_{1,3}}
G_{24}^{(l_{24})}(\gamma_{24})
-{Z_1}^{l_{1,3}-1} H_3^{(l_{124}}(\gamma_{124}) , 
\\
(\gamma_{12}-\gamma_1\gamma_2)
{Z_2}^{l_{234}-1}{(Z_3Z_4)}^{l_{34}-1}
= \gamma_2 (Z_3Z_4)^{l_1-l_{12}}
F_1^{(l_1)}(\gamma_1)
+{Z_1}^{l_{12}}
F_2^{(l_2)}(\gamma_2)
-{Z_2}^{l_{234}-1}
G_{12}^{l_{12}}(\gamma_{12}) .
\end{array}
$$
Hence 
$$
\gamma_{234}=\gamma_{34}\gamma_2 , \
\gamma_2 =\gamma_{23}\gamma_{124} , \ 
\gamma_{124}=\gamma_1\gamma_{24} , \ 
\gamma_{12}=\gamma_1\gamma_2 \ ,
$$
which are the same relations of $u^\prime_I$s in
(\req(Ccip)) for $i=1$ under the identification: 
$u^\prime_1=\gamma_1, u^\prime_2=\gamma_{34}, 
u^\prime_3=\gamma_{24}, u^\prime_4=\gamma_{23}$.
By the similar argument, all the
relations of (\req(Ccip)) hold; therefore 
$J=I(y)$ for
$y$ near
$x_{C_1^\prime}$ having the coordinates
$u^\prime_i$s.  

By the results we have obtained, one concludes
that $\hl^G(\CZ^4)$ is a smooth toric
variety, hence of the form $X_{\Xi^{**}}$ where
$\Xi^{**}$ is a simplicial decomposition 
of $\Delta$ which is refinement of $\Xi^*$
corresponding to the morphism $\lambda$. Indeed,
the above analysis of local structure of ${\rm
Hilb}^G(\CZ^4)$ has shown $\Xi^* = \Xi^{**}$,
therefore $\lambda$ is an isomorphism between 
$X_{\Xi^*}$ and $\hl^G(\CZ^4)$. The
number of exceptional divisors appearing in the
canonical bundle of $X_{\Xi^*}$ is equal to 
$\frac{r(r+1)(r+2)}{6}$
by (\req(numbty)). 
$\Box$ \par \vspace{.2in} \noindent

\section{G-Hilbert Scheme over ${\bf \CZ^3/
{\goth A}_4}$}
It is known that the alternating group ${\goth
A}_{n+1}$ is a simple group except $n=2,3$, in
which cases, ${\goth
A}_3 \simeq \ZZ_3$ and ${\goth
A}_4 $ is isomorphic to the ternary trihedral
group $ (\ZZ_2 \times\ZZ_2) \lhd \ZZ_3$. The
$G$-Hilbert scheme for ${\goth
A}_3$ is the minimal resolution of $\CZ^2/{\goth
A}_3$. In this section we are going to give a constructive proof of the
smooth and explicit crepant structure of
$\hl^{{\goth A}_4}(\CZ^3)$.
\begin{thm}\label{th:A4sm}
$\hl^{{\goth A}_4}(\CZ^3)$ is a
smooth variety with trivial canonical bundle.
\end{thm}
We shall devote the
rest of this section to the proof  of the above
theorem, and always denote $G= {\goth A}_4$. Introduce the following coordinates
$(z_1, z_2, z_3)$ of
$V$ in
$(\req(V))_{n=3}$,
$$
\begin{array}{l}
z_1= - \widetilde{z}_1 + \widetilde{z}_2 +
\widetilde{z}_3 - \widetilde{z}_4 , \ \
z_2= \widetilde{z}_1-\widetilde{z}_2+
\widetilde{z}_3 -\widetilde{z}_4, \ \ 
z_3= \widetilde{z}_1 + \widetilde{z}_2
-\widetilde{z}_3 -\widetilde{z}_4 \ ,
\end{array}
$$
where $\sum_{j=1}^4
\overline{z}_j =0$. The irreducible representation
of $G$ on $\CZ^3 (= V)$, denoted by ${\bf 3}$,
has the following matrix forms for generators
of $G$,
$$
\begin{array}{l}
(12)(34) \mapsto \left( \begin{array}{ccc}
-1 & 0 & 0 \\
0& -1 &0 \\
0& 0 & 1
\end{array} \right)
, \ (13)(24) \mapsto \left( \begin{array}{ccc}
-1& 0 & 0 \\
0& 1&0 \\
0& 0 & -1 
\end{array} \right) , \ \
(123) \mapsto \left( \begin{array}{ccc}
0& 1 & 0 \\
0& 0&1 \\
1& 0 & 0
\end{array} \right) \ .
\end{array} 
$$
There are 4 distinct irreducible $G$-modules,
$
{\rm Irr}(G) = \{ {\bf 1}, {\bf 1}_\omega, {\bf
1}_{\overline{\omega}}, {\bf 3}  \}$, 
where $\omega:= e^{\frac{2 \pi \sqrt{-1}}{3}}$,
and  
$ {\bf
1}_{*}$ is the $G$-character 
determined only by the $(123)$-value $*$. Using the
coordinates $(z_i)_{i=1}^3$ of $\CZ^3$, the generators of 
$G$-invariant polynomials  in $\CZ[Z]$ are: 
$$
\begin{array}{llll}
Y_1 &= Z_1^2 + Z_2^2 + Z_3^2 ,&
Y_2 &= Z_1Z_2Z_3 , \\
Y_3 &= Z_1^2Z_2^2 + Z_2^2Z_3^2 + Z_3^2Z_1^2 ,&
X &= (Z_1^2- Z_2^2)(Z_2^2- Z_3^2)(Z_3^2-Z_1^2) \ .
\end{array}
$$
Note that the above variables are related to
$s_2, s_3, s_4, d$ in $(\req(SA))_{n=3}$ by the relations, $
Y_1= -8 s_2$, $Y_2 = -8 s_3$ , $ Y_3 = 16 s_2^2 -
64 s_4$ , $ X = 64 d $. 
The $G$-invariant polynomial relation (\req(SA)) with $F_3$ in (\req(eqnS)) becomes 
\be
X^2 =   - 4 Y_1^3
Y_2^2 -27 Y_2^4  +18 Y_1Y_2^2Y_3 + Y_1^2
Y_3^2 - 4 Y_3^3 \ .
\ele(eqTri)
Let $\CZ[Z]_j$ be the space of homogeneous
polynomials of degree $j$, and denote $
I(o)^{\bot}_j = I(o)^\bot \cap \CZ[Z]_j$. 
Then $I(o)^{\bot}_j$ is a $G$-submodule of
$I(o)^{\bot}$. In fact, the only non-zero
$I(o)^{\bot}_j$s are among the range $0 \leq j
\leq 5$, whose $G$-irreducible factors are as
follows, ( an equivalent form see, e.g., Table 2.2
in
\cite{GNS}),
\bea(llll)
I(o)^{\bot}_0 =m_0 &  \simeq  &{\bf 1}, & m_0 =
\CZ ,  
\\ I(o)^{\bot}_1 = m_1 &
\simeq  &{\bf 3},  &m_1 = \{ Z_1, Z_2, Z_3
\} ,
\\
I(o)^{\bot}_2 = m_2 +
m_3 +m_4  &
\simeq  &{\bf 1}_{\overline{\omega}} + {\bf
1}_{\omega}+ {\bf 3} ,& m_2 = \{f
\}, \ m_3 = \{ \overline{f} \}, \
m_4 =\{  Z_2Z_3, Z_3Z_1, Z_1Z_2
\} ,  \\
I(o)^{\bot}_3 = m_5 + m_6 &
\simeq  &{\bf 3} + {\bf
3} ,& m_5 = f \{ Z_1 ,
\omega^2 Z_2 , 
\omega Z_3 
\} , m_6= \overline{f} \{ Z_1  ,
\omega Z_2 ,
\omega^2 Z_3  \} ,   \\
I(o)^{\bot}_4 = m_7 +
m_8 + m_9  &
\simeq  &{\bf 1}_{\overline{\omega}} + {\bf
1}_{\omega} + {\bf 3} , & m_7 =\{  
\overline{f}^2 \},  m_8 = \{ f^2  \} , 
m_9= f \{ \omega
Z_1Z_2 ,  Z_2Z_3 , \omega^2 Z_3Z_1\} , \\
I(o)^{\bot}_5 = m_{10} & \simeq & {\bf 3} , &
m_{10} = \overline{f}^2 \{ Z_1 , \omega^2 Z_2 ,
\omega Z_3  \} \ ,
\elea(mon)
where $f:= \sum_{j=1}^3 \omega^{j-1} Z_j^2 $, 
$\overline{f}:= \sum_{j=1}^3 \omega^{2j-2} Z_j^2
$.  We have the
$G$-irreducible decomposition, $I(o)^\bot =
\sum_{k=0}^{10} m_k$.  Note that
$f\overline{f}, f^3, \overline{f}^3$ are
$G$-invariant polynomials with the following
relations,
\bea(l)
f \overline{f} =  Y_1^2 - 3 Y_3 , \ \
f^3 - \overline{f}^3 = 3 (\omega^2 -\omega) X
, \ \
f^3 + \overline{f}^3 = 27 Y_2^2 - 9 Y_1Y_3 + 2
Y_1^3 \ .
\elea(fXY)

\begin{lma}\label{lem:mrel}
Among $m_k$s   $(1 \leq k \leq 10)$, the following 
tree diagram holds:
$$
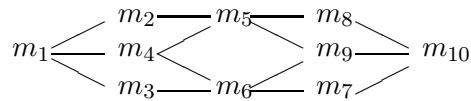
\put(-115,-1){\shortstack{$m_1$}}
\put(-100,2){\line(2,1){20}}
\put(-100,0){\line(1,0){20}}
\put(-100,-2){\line(2,-1){20}}
\put(-75,-1){\shortstack{$m_4$}}
\put(-75,12){\shortstack{$m_2$}}
\put(-75,-15){\shortstack{$m_3$}}
\put(-60, 14){\line(1,0){20}}
\put(-60, -14){\line(1,0){20}}
\put(-60, 0){\line(2,1){20}}
\put(-60,0){\line(2,-1){20}}
\put(-38,12){\shortstack{$m_5$}}
\put(-38,-15){\shortstack{$m_6$}}
\put(-25, 14){\line(1,0){20}}
\put(-25, -14){\line(1,0){20}}
\put(0,-1){\shortstack{$m_9$}}
\put(0,12){\shortstack{$m_8$}}
\put(0,-15){\shortstack{$m_7$}}
\put(-25, 12){\line(2,-1){20}}
\put(-25,-13){\line(2,1){20}}
\put(15,-15){\line(2,1){20}}
\put(15,0){\line(1,0){20}}
\put(15,12){\line(2,-1){20}}
\put(40,-1){\shortstack{$m_{10}$}}
$$
where the $m_j$ of
the  right end of an edge is contained in the
ideal generated by the $m_i$ of the left end
of the edge and $I(o)$.
\end{lma}
{\it Proof.} By the expression of $m_k$, all the
relations in the above diagram are trivial except
the following ones:
\be
m_9 \subset m_6
+ I(o) \ , \ \ m_{10} \subset m_8 + I(o) \ , \ \
m_{10} \subset m_9 + I(o) \ .
\ele(mexc)
Define the irreducible $G$-submodules of
$\CZ[Z]$, isomorphic to ${\bf 3}$: 
$
\overline{m}_{9} = \overline{f} \{ \omega^2
Z_1Z_2 ,
Z_2Z_3 ,
\omega Z_3Z_1  \}$ , $\overline{m}_{10} : =
f^2
\{ Z_1 ,
\omega Z_2 ,
\omega^2 Z_3  \}$. 
Then we have the equalities of ideals in $\CZ[Z]$, $
\langle m_9, I(o) \rangle = \langle \overline{m}_9, I(o) \rangle$, $ 
\langle m_{10}, I(o) \rangle = \langle \overline{m}_{10}, I(o) \rangle$, 
which imply the relations in (\req(mexc)). 
$\Box$ \par \vspace{.2in} \noindent
We shall call
an ideal $J_0$ in
$\hl^G(\CZ^3)$ to be central if $J_0$ is
generated by $I(o)$ and a finite number of
$m_k$s. (The central ideal $J_0$ here will play a similar role of monomial ideals in previous sections for the case of abelian group.) By Lemma \ref{lem:mrel}, there are exactly
four central ideal
$J_0$ with the following $G$-irreducible
decomposition of $\CZ[Z]/J_0$ presented in (\ref{J0A4}).
\bea(ll)
J_0 , & \CZ[Z]/J_0 \\
x_0:= \langle f \rangle +
I(o), 
& m_0 + m_1 +m_3 +m_4 + m_6 +m_7  ; \\
x_0^\prime:= \langle \overline{f} 
\rangle + I(o), & m_0 + m_1 + m_2 +m_4 +m_5 + m_8
; \\
x_\infty := \langle Z_1f, \omega^2 Z_2f,
\omega Z_3f, 
\overline{f}^2  
\rangle + I(o) , & m_0 + m_1 + m_2 +m_3 + m_4
+m_6 ;\\
x_\infty^\prime := \langle Z_1 \overline{f},
\omega Z_2 \overline{f}, \omega^2 Z_3
\overline{f},   f^2
\rangle + I(o) \ , & m_0 + m_1 + m_2 +m_3 + m_4
+m_5 .
\elea(J0A4)
Note that the $J_0$'s presented in (\ref{J0A4}) are
characterized as the ideals  in ${\rm
Hilb}^G(\CZ^3)$ with  monomial polynomial
generators in
$\CZ[Z]$. All the above four elements  lie over
$o \in S_G$ under the morphism
$\sigma_\hl$ of (\req(sHilb)). By the
analysis in \S2.5 of \cite{GNS},  $\sigma_{\rm
Hilb}^{-1}(o)$ consists of a tree of three smooth
rational curves, $L+ l + L^\prime$. Here are the
locations of $J_0$s in
$\sigma_\hl^{-1}(o)$:  $
x_0 \in
(L \setminus l )
\cup L^\prime$, $ x_\infty = L \cap l$, $ 
x_\infty^\prime = L^\prime \cap l$ , $ x_0^\prime \in ( L^\prime
\setminus l ) \cup L $, (see Fig. \ref{fig:A_4}).
\begin{figure}[ht]
$$
\put(-105, 0){\line(3,-1){55}}
\put(-70, -17){\line(3,1){55}}
\put(-28, 0){\line(3,-1){55}}
\put(-105, -4){\shortstack{$\star$ }}
\put(22, -20){\shortstack{$\star$ }}
\put(-65, -25){\shortstack{$x_\infty$ }}
\put(-108, -12){\shortstack{$x_0$ }}
\put(-26, -12){\shortstack{$x_\infty^\prime$ }}
\put(22, -27){\shortstack{$x_0^\prime$ }}
\put(-85, -2){\shortstack{$L$ }}
\put(-45, -5){\shortstack{$l$ }}
\put(0, -5){\shortstack{$L^\prime$ }}
$$
\caption{Tree configuration of  $\sigma_\hl^{-1}(o)$ for $G={\goth A}_4$}
\label{fig:A_4}
\end{figure}
We are going to show that every $J$ in ${\rm
Hilb}^G(\CZ^3)$ can be deformed to one $J_0$ in 
(\req(J0A4)). For $J \in  \hl^G(\CZ^3)$,
denote $h(J)$ be the homogenous ideal in $\CZ[Z]$
generated by the highest total degree part of
elements in $J$. As the top degree of a polynomial
in $\CZ[Z]$ is preserved under the $G$-action,  
$h(J)$ is $G$-invariant. By applying the
Gr\"{o}bner basis technique with a monomial order
of total degree in $\CZ[Z]$, one obtains the same
ideal, ${\rm lt}(J)= {\rm
lt}(h(J))$, hence a set of monomial elements in
$\CZ[Z]$ which represent the basis for both
$\CZ[Z]/J$ and $\CZ[Z]/h(J)$. Therefore $h(J)$
is a homogeneous ideal in $\hl^G(\CZ^3)$.
Note that  $\sigma_{\rm
Hilb}(h(J))=o$. By (2.4) in \cite{GNS}, $h(J)
\in \{x_0, x_0^\prime \} \cup l$. Hence $h(J)$ and
$J$ can be deformed to an  element in
(\req(J0A4)). Now we are going to determine the
local structure near these four central elements
in $\hl^G(\CZ^3)$.

For $J$ near the element $x_\infty$ in ${\rm
Hilb}^G(\CZ^3)$, we have 
$$
J = \langle \overline{f}^2- v_0 f, \ m_5- v_1m_6
-v_2 m_4 - v_3m_1 , \ Y_1 -
\eta_1, Y_2-
\eta_2, Y_3- \eta_3, X - \xi \rangle \ ,
$$
where $(\xi, \eta_1, \eta_2, \eta_3)$ satisfies
(\req(eqTri)), and 
$m_5- v_1 m_6 -v_2 m_4 - v_3 m_1$ is the
$G$-module $\sum_{j=1}^3 \CZ p_j $ with  $
p_1 :=  f Z_1 - v_1 \overline{f} Z_1 - v_2
Z_2Z_3 - v_3 Z_1$, $p_2 :=  f \omega^2 Z_2 - v_1
\overline{f} \omega Z_2 - v_2
Z_3Z_1 - v_3 Z_2 $, $ p_3 :=  f \omega Z_3  - v_1 \overline{f}
\omega^2 Z_3 - v_2
Z_1Z_2 - v_3 Z_3$.
By 
\begin{eqnarray}
f^2 -(v_1 \eta_1 + v_3) \overline{f}
&=& Z_1 p_1 +
\omega^2 Z_2 p_2 +
\omega Z_3 p_3 + v_1(Y_1-\eta_1)
\overline{f}
\in J  \label{relff} \\
(\eta_1 - v_3- v_0 v_1 )f   &=& Z_1
p_1 +
\omega Z_2 p_2 +
\omega^2 Z_3 p_3 + v_1(\overline{f}^2- v_0
f)-(Y_1-\eta_1)f \in J \ , \nonumber
\end{eqnarray}
and the first relation of (\req(fXY)), 
we have
$(3 \eta_3 - \eta_1^2 - v_0 ( v_3 + v_1 \eta_1) )f
\in J $.
As $f
\not\in J$,  we have 
$$
\eta_1 - v_3- v_0 v_1  = 0 \ , \ 
3 \eta_3 - \eta_1^2 = v_0 ( v_3 + v_1 \eta_1) 
\ .
$$
By the relations: $
3Y_2 f - v_2(Z_2^2Z_3^2 + \omega
Z_3^2 Z_1^2 + \omega^2 Z_1^2 Z_2^2) =
Z_2Z_3 p_1 + \omega Z_3Z_1
p_2 +
\omega^2 Z_1Z_2 p_3 \in J$, 
and  $\overline{f}^2 - Y_1 f =
3(Z_2^2Z_3^2 + \omega Z_3^2 Z_1^2 + \omega^2
Z_1^2 Z_2^2)$,   we have 
$$
9\eta_2 + v_2 \eta_1 - v_0 v_2   = 0 \ .
$$
By (\req(fXY)) (\ref{relff}), we have 
$$
3(\omega^2-\omega) \xi \equiv f^3 -
\overline{f}^3 
\equiv (v_1
\eta_1 + v_3 - v_0) f \overline{f} \equiv (v_1
\eta_1 + v_3 - v_0)(\eta_1^2 - 3 \eta_3) \
\pmod{J} \ , 
$$
hence 
$$
3(\omega^2-\omega) \xi = (v_1
\eta_1 + v_3 - v_0)(\eta_1^2 - 3 \eta_3) \ .
$$
By the relation 
$$
\begin{array}{ll}
v_2 p_3 -( \omega- \omega^2)((\omega^2 + \omega
v_1)Z_1p_2 -(\omega +
\omega^2 v_1) Z_2p_1 ) \equiv & \\
v_1(9
\eta_2+ v_0v_1v_2 + v_2v_3 -v_0v_2) Z_3 + ( 3 v_1
\eta_1 + 3 v_3 - 3 v_1v_3 - v_2^2) Z_1Z_2
& \pmod{J}
\ , 
\end{array}
$$
and $Z_1Z_2, Z_3$ representing two
basis elements of $\CZ[Z]/J$, one obtains  
$$
3 v_1
\eta_1 + 3 v_3 - 3 v_1v_3 - v_2^2 = 0 \ , \ 
v_1(9
\eta_2+ v_0 v_1v_2 + v_2v_3 -v_0 v_2) = 0 \ .
$$ 
With all the above relations among $v_j$s,
$\eta_k$s and $\xi$ in the above, one can conclude
that $(v_0, v_1, v_2)$ forms a coordinate system
of
$\hl^G(\CZ^3)$ centered at
$x_\infty$, and the other parameters
in the expression of the ideal $J$ are
expressed by the following relations,
$$
\begin{array}{ll}
v_3  &=  \frac{1}{3} v_2^2 - v_0v_1^2 , \ \ \ \
\ \ \ \ \ \ \  \ \ \ \ 
\eta_1 = \frac{1}{3}v_2^2+v_0v_1 -v_0v_1^2 
,
\\
\eta_2 &= 
\frac{1}{27}v_2(3v_0-3v_0v_1-v_2^2+3v_0v_1^2 ) , 
\ \ \ \
\eta_3 = 
\frac{1}{27}(3v_0v_1^2-v_2^2)
(3v_0v_1^2-3v_0v_1-v_2^2+3v_0) , 
\\
\xi  &=
\frac{\omega-\omega^2}{81}v_0 (v_1+1)
(3v_0v_1^2 +
3v_0-3v_0v_1-v_2^2)(3v_0v_1^3-v_2^2-v_1v_2^2).
\end{array}
$$
Note that  the above
$\xi, \eta_1, \eta_2,
\eta_3$ satisfy the relation (\req(eqTri)). Furthermore, $v_j$s are $G$-invariant
rational functions in $Z_i$s with the following
expressions:
$$
\begin{array}{llll}
v_0 = &
\frac{3(\omega-\omega^2)\xi-9\eta_1\eta_3+
27
\eta_2^2+2\eta_1^3}{2(\eta_1^2-3
\eta_3)} ,
& 
v_1 = &
\frac{(\omega-\omega^2)\xi+ \eta_1\eta_3-
9\eta_2^2}{(\omega-\omega^2)\xi-\eta_1\eta_3+
9\eta_2^2} ,  \\
v_2 = &
\frac{6\eta_2(\eta_1^2-
3\eta_3)}{(\omega-\omega^2)\xi-\eta_1\eta_3+
9\eta_2^2} , &
v_3 = &
\frac{-2\eta_3 (\eta_1^2-
3\eta_3)}{(\omega-\omega^2)\xi-\eta_1\eta_3+
9\eta_2^2} .
\end{array}
$$
This implies $
dZ_1 \wedge dZ_2
\wedge dZ_3 = \frac{\omega-\omega^2}{36} dv_0
\wedge dv_1
\wedge dv_2$.

For $J$ near the element $x_\infty^\prime$ in
$\hl^G(\CZ^3)$, we have 
$$
J = \langle f^2- v_0^\prime \overline{f}, \ m_6-
v_1^\prime m_5 -v_2^\prime m_4 - v_3^\prime m_1 ,
\  Y_1 -
\eta_1, Y_2-
\eta_2, Y_3- \eta_3, X - \xi \rangle \ .
$$
By a similar argument as the case $x_\infty$,
$(v_0^\prime, v_1^\prime, v_2^\prime)$ forms a
coordinate system of
$\hl^G(\CZ^3)$ centered at
$x_\infty^\prime$ with the relations,
$$
\begin{array}{ll}
v_3^\prime & =  \frac{1}{3} v_2^{\prime 2} -
v_0^\prime v_1^{\prime 2} , \ \ \ \ \ \ \ \ \
\ \ \ \ \ \  \ \ \ \
\eta_1 = \frac{1}{3}v_2^{\prime 2}+v_0^\prime
v_1^\prime  -v_0^\prime v_1^{\prime 2}  ,
\\
\eta_2 &= 
\frac{1}{27}v_2(3v_0^\prime
-3v_0^\prime v_1^\prime -v_2^{\prime
2}+3v_0^\prime v_1^{\prime 2} ) , 
\ \ \ \
\eta_3 = 
\frac{1}{27}(3v_0^\prime v_1^{\prime
2}-v_2^{\prime 2}) (3v_0^\prime v_1^{\prime
2}-3v_0^\prime v_1^\prime -v_2^{\prime
2}+3v_0^\prime )
\\
\xi & =
\frac{\omega^2-\omega}{81}v_0^\prime  (v_1^\prime
+1) (3v_0^\prime v_1^{\prime 2}  +
3v_0^\prime
-3v_0^\prime v_1^\prime
-v_2^{\prime 2})(3v_0^\prime
v_1^{\prime 3}-v_2^{\prime 2}-v_1^\prime
v_2^{\prime 2}) \ .
\end{array}
$$
We have 
$$
\begin{array}{llll}
v_0^\prime = &  
\frac{3(\omega^2-\omega)\xi-9\eta_1\eta_3+
27\eta_2^2+2\eta_1^3}{2(\eta_1^2-3
\eta_3)} ,
& 
v_1^\prime = &
\frac{(\omega^2-\omega)\xi+\eta_1\eta_3-
9\eta_2^2}{(\omega^2-\omega)\xi-\eta_1\eta_3+
9\eta_2^2} ,  \\
v_2^\prime = &
\frac{6\eta_2(\eta_1^2-
3\eta_3)}{(\omega^2-\omega)\xi-\eta_1\eta_3+
9\eta_2^2} , &
v_3^\prime = &
\frac{-2\eta_3 (\eta_1^2-
3\eta_3)}{(\omega^2-\omega)\xi-\eta_1\eta_3+
9\eta_2^2} ,
\end{array}
$$
and $
dZ_1 \wedge dZ_2
\wedge dZ_3 = \frac{\omega^2-\omega}{36}
dv_0^\prime \wedge dv_1^\prime
\wedge dv_2^\prime $.

For $J$ near $x_0$ in
$\hl^G(\CZ^3)$, we have 
$$
J = \langle f- u_0 \overline{f}^2, \
\overline{m}_9- u_1 m_6 -u_2 m_4 - u_3 m_1 , \ 
Y_1 -
\eta_1, Y_2-
\eta_2, Y_3- \eta_3, X - \xi \rangle \ ,
$$ 
where $(\xi, \eta_1, \eta_2, \eta_3)$ is as before,
and 
$\overline{m}_9-
u_1 m_6 -u_2 m_4 - u_3 m_1$ is the
$G$-module $\sum_{j=1}^3 \CZ q_j $ with $ 
q_1 :=   \overline{f} Z_2Z_3 - u_1
\overline{f} Z_1 - u_2 Z_2Z_3 - u_3 Z_1 $ , $
q_2 :=  \overline{f} \omega  Z_3Z_1 - u_1
\overline{f} \omega Z_2 - u_2
Z_3Z_1 - u_3 Z_2  $, $
q_3 :=  \overline{f} \omega^2 Z_1Z_2  - u_1
\overline{f}
\omega^2 Z_3 - u_2
Z_1Z_2 - u_3 Z_3 $. 
By the relation, $
-(u_1 + u_0u_3) \overline{f}^2 \equiv - u_1
\overline{f}^2 - u_3 f = Z_1 q_1 +
\omega Z_2 q_2 +
\omega^2 Z_3 q_3 \pmod{J}$, 
we have 
$$
u_1 = - u_0u_3 \ .
$$
By $( 3 \eta_2 - u_1 \eta_1 - u_3) \overline{f}
= Z_1 q_1 + \omega^2 Z_2 q_2 + \omega Z_3 q_3
\in J$, we have 
$$
3 \eta_2 = u_1 \eta_1 + u_3 = u_3( 1- u_0
\eta_1) \ .
$$
By the relations, $f^2 \equiv u_0 \overline{f}^2f \equiv
u_0(\eta_1^2 - 3 \eta_3) \overline{f} \pmod{J}$ and 
$$
\begin{array}{l}
Z_2Z_3 q_1 + \omega^2 Z_3Z_1 q_2 + \omega Z_1Z_2q_3
\equiv (\eta_3 - 3 u_1 \eta_2) \overline{f} -
u_2( Z_2^2Z_3^2 + \omega^2 Z_3^2 Z_1^2 + \omega Z_1^2
Z_2^2 ) \pmod{J} , \\
Z_2^2Z_3^2 + \omega^2 Z_3^2 Z_1^2 + \omega Z_1^2
Z_2^2 \equiv \frac{1}{3}( f^2 - \eta_1
\overline{f}) \equiv \frac{1}{3}( u_0(\eta_1^2 -
3\eta_3) -
\eta_1) 
\overline{f} \pmod{J} \ ,
\end{array}
$$
we have 
$$
(1+u_0u_2) \eta_3 = \frac{1}{3}(9 u_1 \eta_2 -
u_2 \eta_1 + u_0u_2 \eta_1^2) \ .
$$
Using (\req(fXY)), one has 
$$
\begin{array}{l}
u_0 \overline{f}^2f^2 - \overline{f}^3 \equiv
3(\omega^2 - \omega) \xi \ , \ \ \ 
\ \ \ \
2u_0( \eta_1^2 - 3 \eta_3)^2 - 2 \overline{f}^3
\equiv 6(\omega^2 - \omega) \xi  \pmod{J} , \\
2 \overline{f}^3  \equiv 27 \eta_2^3 - 9 \eta_1
\eta_3 + 2 \eta_1^3 - 3(\omega^2 - \omega) \xi
\pmod{J}
\ , 
\end{array}
$$
hence 
$$
3(\omega^2 - \omega) \xi = 2u_0( \eta_1^2 - 3
\eta_3)^2 - 27 \eta_2^3 + 9 \eta_1
\eta_3 - 2 \eta_1^3 \ .
$$
Using the above relations, we have 
$$
\begin{array}{l}
(1+u_0u_2)(1-u_0\eta_1)(Z_1 q_2 - Z_2 q_1 -
u_1(\omega-\omega^2)q_3)\\ 
+ \frac{1}{2+\omega}(1+u_0u_2)\left[
-\omega^2 u_2 Z_3(f- u_0
\overline{f}^2)  +  u_0u_2Z_3( f^2 - u_0(\eta_1
-3\eta_3)\overline{f}\right] \equiv 
(Z_3 \overline{f} -
Z_3 \omega \eta_1) (1-u_0 \eta_1)(
\eta_1+ u_2 + u_0 u_2^2 - 3 u_0^2u_3^2 ) \ 
\pmod{J} \ .
\end{array}
$$
As $Z_3 \overline{f}, Z_3$ are two basis elements
of $\CZ[Z]/J$, their coefficients in the last
term of the above relation are zero. This implies 
$$
\eta_1= - u_2 - u_0 u_2^2 + 3 u_0^2u_3^2 \ .
$$
From all the above relations between $u_i, \eta_j,
\xi$, one concludes that 
$(u_0, u_2, u_3)$ forms a coordinate system
of
$\hl^G(\CZ^3)$ centered at
$x_0$ and  the following relations hold,
$$
\begin{array}{ll}
u_1 & =  - u_0u_3  , \ \ \ \ \ \ \ \ \ \ \ \ \ \
\eta_1 =  - u_2 - u_0 u_2^2 + 3 u_0^2 u_3^2
,
\\
\eta_2 &=  \frac{1}{3}u_3 
(1+ u_0u_2 + u_0^2u_2^2-3 u_0^3u_3^2) ,
\ \ \ \ \
\eta_3 =  \frac{1}{3}(u_2^2-3u_0u_3^2) (1+
u_0u_2 + u_0^2u_2^2 -3 u_0^3u_3^2)
\\
\xi & =
\frac{\omega-\omega^2}{9}(-1+u_0u_2)
(3u_3^2+u_2^3-3u_0u_2u_3^2)(1+ u_0u_2 +
u_0^2u_2^2-3 u_0^3u_3^2) .
\end{array}
$$
Again, the above expressions 
implies the relation (\req(eqTri)), and the
$G$-invariant rational function expression of
$u_i$s are given as follows,
$$
\begin{array}{llll}
u_0 = & 
\frac{6\eta_3- 2\eta_1^2
}{3(\omega^2-\omega)\xi+ 9 \eta_1 \eta_3 - 2
\eta_1^2 - 27 \eta_2^2 } , &  u_1 = &
\frac{\eta_2(-6\eta_3+ 2\eta_1^2)
}{(\omega^2-\omega)
\xi+
\eta_1
\eta_3 - 9
\eta_2^2} , 
\\ u_2 = &
\frac{-(\omega^2-\omega) \xi \eta_1 - 6 \eta_3^2
+ 9 \eta_1 \eta_2^2 +
\eta_1^2 \eta_3 
}{(\omega^2-\omega)
\xi+
\eta_1 \eta_3 - 9
\eta_2^2} , &
u_3 = &
\frac{\eta_2(3(\omega^2- \omega) \xi+
9\eta_1\eta_3 - 27 \eta_2^2 - 2\eta_1^3 )
}{(\omega^2-\omega)
\xi+
\eta_1
\eta_3 - 9
\eta_2^2} ,
\end{array}
$$
hence $
dZ_1 \wedge dZ_2
\wedge dZ_3 = \frac{\omega-\omega^2}{12}
du_0 \wedge du_2
\wedge du_3 $.

For $J$ near the element $x_0^\prime$ in
$\hl^G(\CZ^3)$, we have 
$$
J = \langle \overline{f}- u_0^\prime
f^2, \ m_9- u_1^\prime m_5 -u_2^\prime
m_4 - u_3^\prime m_1 ,
\  Y_1 -
\eta_1, Y_2-
\eta_2, Y_3- \eta_3, X - \xi \rangle \ .
$$
By a similar argument as the case $x_0$, one
obtains that $(u_0^\prime, u_2^\prime,
u_3^\prime)$ is an affine coordinate system with 
$$
\begin{array}{ll}
u_1^\prime & =  - u_0^\prime u_3^\prime  , \ \ \ \ 
\ \ \ \ \ \ \ \
\eta_1 =  - u_2^\prime - u_0^\prime u_2^{\prime
2} + 3 u_0^{\prime 2} u_3^{\prime 2} ,
\\
\eta_2 &=  \frac{1}{3}u_3^\prime 
(1+ u_0^\prime u_2^\prime + u_0^{\prime
2}u_2^{\prime 2}-3 u_0^{\prime 3}u_3^{\prime 2}) ,
\ \ \ \
\eta_3 =  \frac{1}{3}(u_2^{\prime
2}-3u_0^\prime u_3^{\prime 2}) (1+ u_0^\prime
u_2^\prime + u_0^{\prime 2}u_2^{\prime 2} -3
u_0^{\prime 3} u_3^{\prime 2})
\\
\xi & =
\frac{\omega^2-\omega}{9}(-1+u_0^\prime
u_2^\prime ) (3u_3^{\prime 2}
+ u_2^{\prime 3}-3u_0^\prime u_2^\prime
u_3^{\prime 2})(1 + u_0^\prime u_2^\prime +
u_0^{\prime 2} u_2^{\prime 2}-3 u_0^{\prime 3}
u_3^{\prime 2}) , 
\end{array}
$$
and the following relations hold,
$$
\begin{array}{llll}
u_0^\prime = & 
\frac{6\eta_3- 2\eta_1^2
}{3(\omega -\omega^2)\xi+ 9 \eta_1 \eta_3 - 2
\eta_1^2 - 27 \eta_2^2 } , &  u_1 = &
\frac{\eta_2(-6\eta_3+ 2\eta_1^2)
}{(\omega -\omega^2 )
\xi+
\eta_1
\eta_3 - 9
\eta_2^2} , 
\\ 
u_2^\prime = &
\frac{-(\omega -\omega^2 ) \xi \eta_1 - 6 \eta_3^2
+ 9 \eta_1 \eta_2^2 +
\eta_1^2 \eta_3 
}{(\omega -\omega^2)
\xi+
\eta_1 \eta_3 - 9
\eta_2^2} , &
u_3^\prime = &
\frac{\eta_2(3(\omega - \omega^2 ) \xi+
9\eta_1\eta_3 - 27 \eta_2^2-2\eta_1^3 )
}{(\omega -\omega^2 )
\xi+
\eta_1
\eta_3 - 9
\eta_2^2} ,
\end{array}
$$
hence $
dZ_1 \wedge dZ_2
\wedge dZ_3 = \frac{\omega^2-\omega}{12}
du_0^\prime \wedge du_2^\prime
\wedge du_3^\prime$.

With the analysis we have made in this
section, one concludes that $\hl^G(\CZ^3)$
is covered by four 
affine spaces $\CZ^3$ centered at the central
elements in (\req(J0A4)), and the $G$-invariant volume form
$dZ_1\wedge dZ_2 \wedge dZ_3$ of $\CZ^3$
induces a never-vanishing global volume
form of
$\hl^G(\CZ^3)$. This completes the proof
of Theorem \ref{th:A4sm}.

\section{Concluding Remarks} 
In this article, we have provided a detailed
derivation of the smooth toric structure of ${\rm
Hilb}^{A_r(4)}(\CZ^4)$. Its relation with
crepant resolutions of $\CZ^4/A_r(4)$ has been
found, and   different crepant
resolutions connected by flops of 4-folds
can be visualized in the process. We have also
given a constructive verification of the smooth
and crepant properties of ${\rm
Hilb}^{\goth{A}_4}(\CZ^3)$ by a direct
computation method. In the abelian case
$A_r(4)$, the solution has been given in Sects. 3,
4 by the standard toric method, a combinatorial
mechanism built upon monomials in
$\CZ[Z]$, which can be regarded as characters
of the whole torus group $T_0$, containing
$A_r(4)$ as a finite subgroup. The smooth toric
structure of
$\hl^{A_r(4)}(\CZ^4)$ is derived from a
procedure, which mainly consists of two 
steps: first, one obtains a complete list of
monomials ideals in
$\hl^{A_r(4)}(\CZ^4)$ which correspond to
the 0-dimensional toric orbits,( see (\req(A2))
(\req(J0)) (\req(lsrel))); second, by the Gr\"{o}bner
basis technique and a detailed analysis of the
$G$-regular module property of $\CZ[Z]/J$ for an
ideal $J$ in $\hl^{A_r(4)}(\CZ^4)$, one
proceeds to identify the toric coordinates from
the ideal-generators of $J$. In this
manner, the explicit form of the canonical bundle
of 
$\hl^{A_r(4)}(\CZ^4)$  can be determined as
a  disjoint sum 
of exceptional divisors, each of which  could be
blown down to give rise to crepant resolutions
of $\CZ^4/A_r(4)$. These
crepant resolutions are connected by a sequence
of flops in 4-folds through ${\rm
Hilb}^{A_r(4)}(\CZ^4)$. We intend to apply a
similar mechanism to the non-abelian case $G=
{\goth A}_{n+1}$, but relying only
on the data of $G$-representations in $\CZ[Z]$, 
a ``big" group like the torus in the abelian case
does not exist in the  latter case though. In
\S6, we have made a detailed study on the
structure of 
$\hl^{{\goth A}_4}(\CZ^3)$,
which would serve us as a demonstration of the
effectiveness of the method even though its
crepant smooth conclusion is known by now
\cite{BKR}. We have succeeded to give an
explicit verification of  the
crepant smooth structure of 
$\hl^{{\goth A}_4}(\CZ^3)$ following
our  thought by a direct constructive method
via group representations. A similar
analysis to the higher dimensional cases is now
under progress and  partial results are
promising. As to the role of
$G$-Hilbert scheme in the study of crepant
resolution of $S_G$, our conclusion for the case
$G= A_r(4)$ has indicated the non-crepant
property of
$\hl^G(\CZ^4)$, but with a intimate
relation with crepant resolutions of $S_G$. For
higher dimensional case, this kind of link
between
$\hl^G(\CZ^n)$ and some possibly existing
crepant resolutions of $S_G$ could be 
further loosely related. However, the
$G$-Hilbert scheme would still be worth for
further study on its own right due to the
built-in character 
of group representations into the geometry of
orbifolds. This could be a promising
direction of the geometrical study of singularity.
Such program is now under our consideration for
the future study.

\section*{ Acknowledgments}
This work was reported by the second author in the workshop ``Modular Invariance, ADE, Subfactors and Geometry of Moduli Spaces", Kyoto, Japan, Nov. 25- Dec. 2, 2000, and as part of the subject of an Invited Lecture at ``7th International Symposium on Complex Geometry", Sugadaira, Japan, Oct. 23- 26, 2001, for which he would like thank them for their invitation and hospitality.  The research of this paper is supported in part by the National Science Council of Taiwan under grant No.89-2115-M-001-037.


\begin{thebibliography}{99}
\bibitem{BM} J. Bertin and D. Markushevich,
Singlarit\'{e}s quotients non ab\'{e}liennes de
dimension 3 et vari\'{e}t\'{e} de Calabi-Yau,[
Three-dimensional nonabelian quotient singularities and 
Calabi-Yau manifolds]
Math. Ann. 299 (1994) 105-116.  


\bibitem{BKR} T. Bridgeland, A. King and M. Reid, 
The McKay correspondence as an equivalence of derived categories, J. Amer. Math. Soc. {\bf 14}(2001), no. 3, 535-554.

\bibitem{CR} L. Chiang and S. S. Roan, Orbifolds
and finite group representations,Internat. J. of Math. and Math. Sci. 26:11(2001)649-669; math.AG/0007072.



\bibitem{CLO} D. Cox, J. Little and D. O'Shea,
{\it Ideals, varieties, and Algorithms,}
Springer-Verlag, New York-Berlin-Heidelburg, 1992.


\bibitem{D} V. I. Danilov, The geometry of toric varieties, Uspekhi Mat. Nauk 33 (1978), no. 2(200), 85-134.(Russ. Math. Surveys 33:2 (1978) 97-154.)


\bibitem{GNS} Y. Gomi, I. Nakamura and K. Shinoda, Hilbert schemes of $G$-orbits in dimension three,  Asian J. Math. 4 (2000) 51-70.



\bibitem{INj} Y. Ito and H. Nakajima, McKay
correspondence and Hilbert schemes in dimension
three, Topology 39 (2000) 1155-1191, math.AG/9803120.


\bibitem{IN1} Y. Ito and I. Nakamura, McKay
correspondence and Hilbert schemes, Proc.
Japan Acad. 72 (1996) 135-138

\bibitem{IN2} Y. Ito and I. Nakamura, Hilbert  schemes
and simple  singularities, New trends in
Algebraic Geometry (Proc. of the July 1996
Warwick European Alg. Geom. Conf.), London Math. Soc. Lecture Note Ser., vol. 264, Cambridge Univ. Press (1999) 151-233.


\bibitem{M} G. Kempf, F. Knudson, D. Mumford
and B. Saint-Donat, {\it Toroidal embedding} 1 ,
Lecture Notes in Math., Vol. 339,
Springer-Verlag, New York, 1973.




\bibitem{Kl} F. Klein, {\it Gesammelte
Mathematische Abhandlungen.},
Springer-Verlag 1922 (reprint 1973).


\bibitem{MBD} G. A. Miller, H. F. Blichfeldt and L. E.
Dickson, {\it Theory
and applications of finite groups}, John Wiley and
Son, New York, 1916.


\bibitem{Mo} S. Mori, Birational classification of algebraic
threefolds, Proceedings of the Internationational Congress of Mathematicians, Vol. I, II (Kyoto 1990), Math. Soc. Japan, Tokyo(1991)235-248.


\bibitem{N1} I. Nakamura, Hilbert scheme and simple
singularities
$E_6, E_7$ and $E_8$, Hokkaido Univ. Preprint Series in Math.
No. 362 (1996), Hokaido University, Japan, 1996.


\bibitem{N} I Nakamura, Hilbert schemes of
abelian group orbits, J. Alg. Geom. 10(4)(2001) 757--779.

\bibitem{O} T. Oda, {\it Torus embeddings 
and applications,} Tata Institute of Fundamental Research Lectures on  Mathematics and Physics, vol. 57, Tata Institute of Fundamental Reseach, Bombay(1978).


\bibitem{R89} S. S. Roan, On the generalization 
of Kummer surfaces, 30
(1989) 523-537.


\bibitem{Rtop} S. S. Roan,  Minimal resolutions
of Gorenstein orbifolds in dimension three,
Topology 35 (1996) 489-508.


\bibitem{R97} S. S. Roan, Crepant resolution and
fibred CY manifolds, Preprint MIAS 97-1, Inst. of
Math. Acad. Sinica, Taiwan, 1997.


\end{thebibliography}
\end{document}